\documentclass[amssymb,amsmath,nofootinbib]{revtex4}
\usepackage{amssymb,graphicx}
\usepackage[dvips]{color}

\newcommand{\ds}{\displaystyle}
\newcommand{\om}{\omega}
\newcommand{\omp}{{\textstyle\omega'}}
\newcommand{\ompp}{{\textstyle\omega''}}
\newcommand{\Kappa}{\symbol{26}}
\newcommand{\mfrac}[2] 
{\raisebox{0.045em}{\mbox{\footnotesize$\displaystyle
\frac{#1}{#2}$}}}

\begin{document}

\title{On uniformization of Burnside's curve $y^2=x^5-x$}

\author{Yu.\,V.\,Brezhnev}


\email{brezhnev@mail.ru}




\begin{abstract}
{\footnotesize Main objects of uniformization of the curve
$y^2=x^5-x$ are studied: its Burnside's parametrization,
corresponding Schwarz's equation, and accessory parameters. As a
result we obtain the first examples of solvable Fuchsian  equations
on torus and exhibit number-theoretic integer $q$-series for
uniformizing functions, relevant modular forms, and analytic series
for holomorphic Abelian integrals. A conjecture of Whittaker for
hyperelliptic  curves and its hypergeometric reducibility are
discussed. We also consider the conversion between Burnside's and
Whittaker's uniformizations.}
\end{abstract}

\noindent {\footnotesize{\sc Journ.\;Math.\;Phys.} {\bf50}(10),
00--00 (2009)\hfill {\tt http://arXiv.org/math.CA/0111150}}

\bigskip

\maketitle

\section{Introduction}
Presently, one example of an algebraic curve of genus $g>1$
\begin{equation}\label{bs}
y^2=x^5-x
\end{equation}
is known where one is able to present all the key objects of its
uniformization in an explicit form. These include: 1)
parametrization $x(\tau),\,y(\tau)$ in terms of known special
functions; 2)  closed differential calculus associated with these
functions; 3) pictures of conformal representations resulting from
the uniformizing functions; 4) differential equations on functions
$x(\tau)$ and $y(\tau)$,  their solutions, and; 5) matrix
representations for monodromy groups of related Fuchsian equations.

Somewhat surprising facet is the fact that, except for
parametri\-za\-tions by means of modular functions (the theory of
modular equations \cite{klein,erdei}), the straightforward statement
of these questions was not considered in the literature, including
even the famous Klein curve $y^7=x^3-x^2$ \cite{eightfold}. The
parametrization property for the curve (1) was found by William
S.\,Burnside  in 1893 \cite{burnside},  but this result has,
however, received almost no mention in the papers and has not
appeared in the monographic literature devoted to uniformization,
automorphic functions or other relevant material. To the best of our
knowledge, since 1893 the work \cite{burnside} has been  mentioned
by R.\;Rankin in 1958 (see e.\,g. \cite{rankin2}) and comparatively
recently he drew attention to it in work \cite{rankin}. On the other
hand, explicit instances of all of these objects in the case of
higher genera would enable one to construct nice applications
wherever algebraic functions and Riemann surfaces appear. Such
topics include algebraic-geometric integration, conformal field
theories \cite{tah}, integrable quantum/classical dynamical systems
and nonlinear {\sc pde}'s, equations of Picard--Fuchs
\cite{ford,mckay2}, second order linear ordinary differential
equations ({\sc ode}'s) of Fuchsian type \cite{mckay1}, number
theory, ``Monstrous Moonshine'' \cite{gannon}, and many others. Yet
another and  deep application of uniformization originated in works
of Takhtajan and collaborators in the late 1980's and for further
applications of the general uniformization theory see work
\cite{tah} and references therein. However, it is pertinent to note
that the general theory has long experienced the lack of nontrivial
illustrative instances.

History of the curve (1) goes back to the work of Bolza \cite{bolza}
where he obtained period matrix for this curve and  its automorphism
group. Geometry of fundamental domains of Fuchsian groups
(hyperbolic polygons)  for curves of lower genera, including the
curve (1), has been studied extensively in the literature and rather
well developed. See for example works \cite{kuusalo,eightfold} and
references therein. In 1958 Rankin found the correct Fuchsian
equation uniformizing this curve in the framework of an approach of
Whittaker \cite{whittaker1}. More recently, Rankin simplified
\cite{rankin} (inserting radicals however) Burnside's parametrizing
functions (\ref{burnside}) and considered some of their group
properties: that is,  the structure of automorphism groups of the
functions $x(\tau)$ and $y(\tau)$. The group of the function
$x(\tau)$ was described earlier by Klein \& Fricke \cite{klein} and
used by Rankin  and Burnside himself \cite{burnside}. At present,
this is as far as we known about uniformization of the curve (1).

In this note, we would like to draw attention to this remarkable
example in a classical framework\footnote{See for example last
sentence in the book \cite{ford} and Weyl's emphasis on
pp.\,176--177 in the very first monograph on Riemann surfaces
\cite{weyl}.} and to exhibit the differential properties of
Burnside's Riemann surface. It is widely known that these properties
are fully determined  by the second  order linear differential
equations of Fuchsian class on plane \cite{ford} (Secs.\,II, III) in
so far as the monodromical groups of these equations give, broadly
speaking, a matrix $2\!\times\!2$-representation of fundamental
group $\pi_1^{}$ of the Riemann surface. If the surface (a curve)
may be realized as a cover of a Riemann surface of genus $g=1$
(torus) then these equations can be transformed to Fuchsian
equations on this torus. We explain this in Sec.\,IV and consider
also solutions of these equations. In addition to this, in
Sec.\,III, we discuss relation with a conjecture of Whittaker
\cite{whittaker2,whittaker3}. Automorphic forms and analytical
integer $q$-series for functions we construct are considered in
Sec.\,V. All the results are new.

As well as being the first completely describable example,
Burnside's curve also provides rich material for different kind
nontrivial generalizations and observations, some of which we
expound below. Throughout the paper we have adhered to basic and
standard terminology in the theory of automorphic functions and
uniformization \cite{ford,erdei,weyl,mckay1}. The classical
bibliography is listed in the book \cite{ford} and  the works
\cite{mckay1,mckay2,br2} provide additional and modern references.

Let $x=x(\tau)$ be  meromorphic  automorphic function on a Riemann
surface of the algebraic curve $F(x,\,y)=0$ with a local
uniformizing parameter $\tau$.  Then $x(\tau)$ satisfies the
following nonlinear autonomous differential equation of the third
order \cite{ford,poincare,rankin2} (we call it Schwarz's equation)
\begin{equation}\label{Ds1}
\frac{\big\{x,\,\tau \big\}}{x_\tau^2}=Q(x,y)\,, \qquad
\big\{x,\,\tau\big\}:= \frac{x_{\tau\tau\tau}}{x_{\tau}}-\frac32\!
\left(\! \frac{x_{\tau\tau}}{x_{\tau}}\!\right)^{\!\!2},
\end{equation}
wherein $Q(x,y)$ is some rational function of $x$ and $y$
\cite{ford,poincare}. To determine its coefficients, so that
monodromy group of the associated Fuchsian equation
$$
\Psi_{\mathit{xx}}=\frac12\,Q(x,y)\,\Psi\,,\qquad\tau
=\frac{\Psi_1(x)}{\Psi_2(x)}
$$
shall be Fuchsian, is the celebrated problem of accessory parameters
considered in the context of algebraic curves
\cite[pp.\,222--228]{poincare}.  With the correct parameters in hand
the quantity $\tau$ becomes the global uniformizing parameter.

Compared to the theory elliptic functions,  currently available
analytical description  of uniformizing functions for  the genera
higher than unity is poorly developed. By this we mean: 1)
determining {\sc ode}'s (\ref{Ds1}), their solutions, and correct
accessory parameters; 2) effective series expansions and numerical
computations; 3) inversion problems in fundamental polygons, i.\,e.
search for solutions $\tau$ to equations of the type $x(\tau)=A$; 4)
conformal representations; 5) Abelian integrals as functions of the
global parameter $\tau$; 6) addition theorems for these Abelian
integrals and relation with the Jacobi inversion problem.  In the
following sections we shall fill up some of these gaps in the
example of the curve (1). The availability of all of these
attributes in the case of genus $g=1$ provides the great efficiency
of  elliptic functions and their numerous applications, which cannot
be said of the cases $g>1$ for reasons of the insufficiently
advanced analytic tools.

In this paper we follow classical Burnside's $\wp$-formulae,
although they have a simpler representation in terms of jacobian
$\vartheta$-constants. Such and exhaustive $\vartheta$-function
description of the example (1) (without intersections with the
present work)  have been detailed in the work \cite{br2}. In the
same place further generalizations and  extended bibliography are
presented.

\section{Burnside's uniformization and the Schwarz equation}

 In this section we show that Burnside's result provides
the first example of nontrivial algebraic curve where one is known
both corresponding Fuchsian equation with correct accessory
parameters and   uniformizing functions, and also complete
differential calculus of these functions. We reproduce a
parametrization of the curve (\ref{bs}) in an explicit form more
suitable for our purpose. Let $\tau$ be a complex variable with
nonzero imaginary part $\boldsymbol{\Im}(\tau) \ne 0$ and
Weierstrass's $\sigma,\,\zeta,\,\wp,\,\wp'$-functions are taken with
half-periods $(\omega,\,\omp)=(2,\,2\tau)$.  For example
$$
\wp(z|2,2\tau)= \wp\big(z;\,g_2^{}(2,2\tau),\,g_3^{}(2,2\tau)\big)
\qquad \big(=: \wp(z)\big)\,.
$$
In this notation the parametrization of Burnside \cite{burnside} has
the form
\begin{equation}\label{burnside}
\left\{\!
\begin{array}{l} \ds
x= \frac{\wp(1)-\wp(2)}{\wp(\tau)-\wp(2)}\\
\\
\ds y=4i\,\frac {\big[\wp(\tau)\!-\!\wp(2\tau)\big]\!
\big[\wp\big(\frac{\tau}{2}\big)\!-\!\wp(\tau)\big]\!
\big[\wp\big(\frac{\tau}{2}\big)\! -\!\wp(\tau\!+\!2)\big]\!
\big[\wp\big(\frac12\big)\!-\!\wp(2\tau\!+\!1)\big]\!
\big[\wp\big(\frac12\big)\!-\!\wp(1)\big]_{{\mathstrut}}}
{\big[\wp\big(\frac \tau 2\big)\!-\!\wp(1)\big]^{{\mathstrut}}\!
\big[\wp\big(\frac \tau 2\big)\!-\!\wp(2\tau\!+\!1)\big]
\wp'\big(\frac 12\big)\,\wp'(\tau)}
\end{array}\right.\!.
\end{equation}

{\bf Proposition\;1.} {\em The Schwarz equation {\rm (\ref{Ds1})}
and, therefore, the accessory parameters for  Burnside's
parametrization {\rm (\ref{bs})}, {\rm(\ref{burnside})} have  the
form}
\begin{equation}\label{Ds2}
\frac{\big\{x,\,\tau\big\}}{x_\tau^2}=  -\frac12\left\{
\frac{1}{x^2}+\frac{1}{(x-1)^2} +\frac{1}{(x+1)^2}+
\frac{1}{(x-i)^2}+\frac{1}{(x+i)^2}
-\frac{4\,x^3+0}{x^5-x}\right\}\,.
\end{equation}
Of interest is to give a direct proof, although the verification of
the calculations is not straightforward because the
$\zeta,\,\wp,\,\wp'$-functions are not closed under differentiation
with respect to $\tau$. For example, we shall require a closed
system of differential equations satisfied by
$g_{2,3}^{}(\om,\,\omp)$ and the periods $\eta,\eta'(\omega,\,\omp)$
of the elliptic integral $\zeta(z|\om,\omp)$:
\begin{equation}\label{diffg2g3}
\frac{dg_2^{}}{d\tau}= \frac{i}{\pi}
\Big(8\,g_2^{}\,\eta-12\,g_3^{}\Big),\qquad \frac{d g_3^{}}{d\tau}=
\frac{i}{\pi} \Big(12\,g_3^{}\,\eta-\mfrac23\,g_2^2\Big),\qquad
\frac{d\eta}{d\tau}=
\frac{i}{\pi}\Big(2\,\eta^2-\mfrac16\,g_2^{}\Big),
\end{equation}
where we use the notation $g_{2,3}^{}=g_{2,3}^{}(\tau)$ for
$g_{2,3}^{}(1,\,\tau)$ and the same for $\eta$ and $\eta'$. Some of
the intermediate results in the proof will be required later. Since
the arguments of Weierstrass's functions will constantly appear in
the following calculations, we adopt the concise notation
$\wp_\alpha^{} := \wp(\alpha|2,2\tau)$ etc. Throughout the paper in
the calculations that follows (being automated with a computer), we
use intensively
known and new
properties of elliptic and modular
functions collected in a form of reference source in work
\cite{br3}.

{\em Proof.\/} A straightforward substitution of (\ref{burnside})
into (\ref{Ds1}) generates, apart from cumbrous formulae, many
variables, including $\tau$ in an explicit way. But they are not
algebraically independent. Using the addition theorems for the
functions $\zeta(2\alpha)$, $\wp(2\alpha)$, $\wp'(2\alpha)$ at
points $\alpha=1$ and $\tau$, we obtain six identities.  We shall
need also important formulae for the invariants $g_{2,3}^{}$
$$
g_2^{}=4\,\frac{\wp_\alpha^3-\wp_\nu^3}{\wp_\alpha^{}-\wp_\nu^{}}-
\frac{{\wp'_\alpha}^{\!2}-{\wp'_\nu}^{\!2}}{\wp_\alpha^{}-\wp_\nu^{}}\,,
\qquad g_3^{}=\frac{\wp_\nu^{}\,{\wp'_\alpha}^{\!2}-
\wp_\alpha^{}\,{\wp'_\nu}^{\!2}}{\wp_\alpha^{}-\wp_\nu^{}}-
4\,\wp_\alpha^{}\,\wp_\nu^{}\,(\wp_\alpha^{}+\wp_\nu^{})\,,
$$
derivable  from the relations
${\wp'_{\alpha,\nu}}^{\!\!\!\!\!\!2}=4\,\wp_{\alpha,\nu}^3-g_2^{}\,
\wp_{\alpha,\nu}^{}-g_3^{}$, and which are valid under $\alpha \ne
\nu$.  We use them with $\alpha=1$ and $\nu=2$.  Now two of the six
above identities are transformed to $0=0$. The four remaining, after
some simplification, turn into the following identities
\begin{equation}\label{four}
\begin{array}{c}
(\eta-4\,\zeta_1^{})^2=8\,\wp_1^{}+4\,
\wp_2^{}{}_{{\ds\mathstrut}}\,,\qquad \wp'_1=(\eta-4\,\zeta_1^{})\,
(\wp_1^{}- \wp_2^{{\ds\mathstrut}})
_{{\ds\mathstrut}}\,, \\
\wp'_\tau=\ds \frac{2^{\ds\mathstrut}}{\eta'-4\,
\zeta_{\tau_{\ds\mathstrut}}^{}}\,
\big(3\,\wp_\tau^2+\wp_1^2-2\,\wp_1^{}\,\wp_2^{}-2\,\wp_2^2\big)\,, \\
\wp_\tau^4+2\,\wp_2^{}\, \wp_\tau^3 +
6\big(\wp_1^2-2\,\wp_1^{}\,\wp_2^{}-\wp_2^2\big)\wp_\tau^2 -
2\,\wp_2^{}\big(3\,\wp_1^2-6\,\wp_1^{}\,\wp_2^{}-4\,\wp_2^2\big)
\wp_{\tau{}_{\ds\mathstrut}}^{{}^{\mathstrut}}+{}\\
{}+\wp_{1}^{4^{\ds\mathstrut}}-4\,\wp_1^3\,\wp_2^{} +6\,\wp_1^2\,
\wp_2^2-
4\,\wp_1^{}\,\wp_2^3-4\,\wp_2^4=0\,.\qquad\qquad\qquad\qquad\quad\;\,
\end{array}
\end{equation}
The equations (\ref{diffg2g3})--(\ref{four}), together with the
rules of differentiation of  Weierstrassian functions \cite{br3},
contain all the information required for the proof. Using the first
formula in (\ref{burnside}) we obtain the following useful
expressions:
\begin{equation}\label{p2}
\wp'_{\tau}=\frac{144}{4\,\zeta_\tau^{}-\eta'}\,
\frac{x^4-1}{(x^4-6\,x^3+6\,x^2-6\,x+1)^2}\,, \qquad
\wp_1^{}=\frac{x^4-6\,x^3+6\,x^2-6\,x+1}{x^4+6\,x^2+1}\,\wp_2^{}\,,
\end{equation}
\smallskip
\begin{equation}\label{g2}
g_2^{}(\tau)= 2^6
\,3\,\frac{x^8+14\,x^4+1}{(x^4+6\,x^2+1)^2}\,\wp_2^2\,,\qquad\quad
g_3^{}(\tau)=
-2^9\,\frac{x^{12}-33\,(x^8+x^4)+1}{(x^4+6\,x^2+1)^3}\,\wp_2^3\,,
\end{equation}
whereupon three derivatives of the $x(\tau)$-function acquire the
form
\begin{eqnarray}
\label{xt} x_{\tau}&=&\ds\phantom{-}
\frac{24}{\pi}i\,\frac{x^5-x}{x^4+6\,x^2+1}\,\wp_2^{}\,,\\ \nonumber \\
x_{\tau\tau}&= &\ds\nonumber -\frac{96}{\pi^2}\,
\frac{(x^4+6\,x^2+1)\,\eta+2\,(5\,x^4-1)\,\wp_2^{}}
{(x^4+6\,x^2+1)^2}\,(x^5-x)\,\wp_2^{}\,,\\ \nonumber \\
x_{\tau\tau\tau}& = & \ds\nonumber -\frac{576}{\pi^3}i\left\{ \frac
{(x^4+6x^2+1)\,\eta+4\,(5x^4-1)\,\wp_2^{}} {(x^4+6x^2+1 )^2}\,\eta
+8\frac{(11x^8-26x^4-1)}{(x^4+6x^2+1)^3}\,\wp_2^2
\right\}\!(x^5\!-\!x)\,\wp_2^{}.
\end{eqnarray}
The formula (\ref{Ds2}) is obtained after substitution of these
derivatives into (\ref{Ds1}). \hfill $\blacksquare$

{\em Remark\;$1$\/}. We can view the last identity in (\ref{four})
as a plain curve in projective coordinates
$\big(\wp_\tau^{}\!:\!\wp_1^{}\!:\!\wp_2^{} \big)$. It has genus
$g=0$. First formula in (\ref{burnside}) and second one in
(\ref{p2}) yield
\begin{equation}\label{pt}
\frac{\wp_\tau^{}}{\wp_2^{}}=\frac{x^4-5}{x^4+6\,x^2+1}\,, \qquad
\frac{\wp_1^{}}{\wp_2^{}}=\frac{x^4-6\,x^3+6\,x^2-6\,x+1}{x^4+6\,x^2+1}
\end{equation}
and therefore the quantity $x$ may be considered as a global
parameter in rational parametrization of this curve. Variation of
$x\in\overline{\mathbb{C}}$ is equivalent, through (\ref{burnside}),
to variation of $\tau\in\mathbb{H}^+$ in fundamental polygon for
group $\boldsymbol{\Gamma}(4)$ \cite{rankin,br2} and thereby correct
variation of quantities $(\wp_\tau^{},\wp_1^{},\wp_2^{})$.

\section{Fuchsian equations and a conjecture of Whittaker}
\subsection{Fuchsian equations associated to Burnside's parametrization}
Integrability of the Fuchsian equation associated with the formula
(\ref{Ds2})
\begin{equation}\label{last}
\begin{array}{l}
\Psi_{\mathit{xx}}=\ds
\;\;\,\frac12\,Q(x)\, \Psi={}\\ \\
\phantom{\Psi_{\mathit{xx}}}=\ds
-\frac14\,\frac{x^8+14\,x^4+1}{(x^5-x)^2}\,\Psi
\end{array}
\end{equation}
immediately follows from known properties of the Schwarzian
(\ref{Ds1}). Denoting for brevity $\frac{1}{\!\sqrt{z}\,}$ as
$\sqrt[\leftroot{-1}\uproot{2}-2]{z}$ we have the well-known
identity \cite{ford}
$$
\Big(\!\! {\textstyle\sqrt[\leftroot{-2}\uproot{2}-2]{\tau_{\!x}^
{\phantom{i}}}}\Big)_ {\!\mathit{xx}}=\frac12\,Q(x) \cdot\!\!\!
{\textstyle\sqrt[\leftroot{-1}\uproot{2}-2]
{\tau_{\!x}^{\phantom{i}}}}\,.
$$
Thus, setting $\Psi(x)=
\!\!\sqrt[\leftroot{-1}\uproot{2}-2]{\tau_{\!x}^{\phantom{i}}}$ we
obtain an integral of equation (\ref{last}). There are numerous
forms for integrals of Fuchsian solvable equations and we, following
(\ref{burnside}), give them in Burnside's $\wp$-manner.

{\bf Proposition\;2.} {\em The general multi-valued integral of
equation {\rm(\ref{last})} is given by the formula
\begin{equation}\label{x45}
\Psi(x)=
\sqrt{\mfrac{x^5-x}{x^4-5}}\,\sqrt{\wp\big(\tau(x)|2,2\tau(x)\big)}
\,\big(A\,\tau(x)+B \big)\,,
\end{equation}
where  function $\tau(x)$ is determined by the inversion of the
expression}
\begin{equation}\label{xt12}
x=\frac{\wp(1|2,2\tau)-\wp(2|2,2\tau)}
{\wp(\tau|2,2\tau)-\wp(2|2,2\tau)}\,.
\end{equation}
Proof follows from the formulae (\ref{xt}) and (\ref{pt}). The
second Fuchsian equation is an equation on the function
$\widetilde\Psi(y)$ defined by the rule $\widetilde\Psi(y)=\!\!
\sqrt[\leftroot{-2}\uproot{2}-2]{\phantom{\tau_{\!x}^{\phantom{i}}}}
\!\!\!\!\!\tau_{\!y}\,$. It has a form of equation with algebraic
coefficients:
\begin{equation}\label{Qy}
\begin{array}{l}
\widetilde\Psi_{\mathit{yy}}=\ds
\;\;\,\frac12\,Q(x,y)\, \widetilde\Psi={}\\ \\
\phantom{\Psi_{\mathit{yy}}}=\ds
-\frac14\,\frac{5^4\,x\,y^6+415\,x^2\,y^4-511\,x^3\,y^2+
255\,x^4+1}{(5^4\,x\,y^6+1375\,x^2\,y^4+1025\,x^3\,y^2+
255\,x^4+1)\,y^2}\,\widetilde\Psi\,.
\end{array}
\end{equation}
General representation for its integral is given by the formula
$$
\widetilde\Psi(y)=
{\textstyle\sqrt{y_\tau^{\phantom{a}}}}\,\big(A\,\tau(y)+B\big)\,,
$$
where $\tau(y)$ is an inversion of the second formula in
(\ref{burnside}).  Schwarz's equation corresponding to the second
function $y(\tau)$ has the form (\ref{Ds1}) with  $Q(x,y)$-function
defined by the right hand side of  expressions (\ref{Qy}).

The group of automorphisms $\boldsymbol{\mathrm{Aut}}(x(\tau))$ of
Burnside's function $x(\tau)$ is of index 24 subgroup in the full
modular group $\mathrm{PSL}_2(\mathbb{Z})$ having the Klein
invariant $J(\tau)$ as a Hauptmodulus \cite{klein,burnside}. Indeed,
with the proof of Proposition\;1 we obtain the formula
\begin{equation}\label{J}
\frac{g_2^3}{g_2^3-27\,g_3^2}=\frac
{1}{108}\,\frac{(x^8+14\,x^4+1)^3}{(x^5-x)^4}=J(\tau)\,.
\end{equation}
The factor group
$\mathrm{PSL}_2(\mathbb{Z})/\boldsymbol{\mathrm{Aut}}(x(\tau))$ is
the octahedral group \cite{klein} of order 24 \cite{bolza}. Taking
into account a permutation of sheets $y(\tau+4)=-y(\tau)$
\cite{rankin}, it becomes the maximal automorphism group of order 48
for the curves of genus two.

\subsection{Whittaker's conjecture}
In 1929 E.\,Whittaker \cite{whittaker2} proposed a pure algebraic
solutions of the transcendental problem of accessory parameters for
certain Fuchsian groups of the first kind without parabolic edges
\cite{ford}. Namely, from words of his son J.\,M.\,Whittaker
\cite{whittaker3}, he suggested that for hyperelliptic algebraic
curves
\begin{equation}\label{hyp}
y^2=(x-e_1^{}) \cdots (x-e_{2g+1}^{})\quad \big(=:f(x)\big)
\end{equation}
the $Q$-function is given by the  formula
\begin{equation}\label{conj}
Q(x)=-\frac38\left\{\frac{f_x^2}{f^2}-\frac{2\,g+2}{2\,g+1}\,
\frac{f_{\mathit xx}^{}}
{f}\right\}.
\end{equation}
The conjecture was checked for some curves in the 1930's by
Whittaker's collaborators \cite{dhar,whittaker3}. What is more, by
considering the hyperlemniscate algebraic curve
\begin{equation}\label{whit}
y^2=x^5+1\,,
\end{equation}
Whittaker \cite{whittaker2} reduced  the associated Fuchsian linear
ordinary differential equation
\begin{equation}\label{fuchs}
\Psi_{\mathit{xx}}= -\frac{3}{16}\left\{
\sum\limits_{k=1}^{5}\frac{1}{(x-e_k^{})^2}-
\frac{4\,x^3+A_1\,x^2+A_2\,x+A_3}{(x-e_1^{})\cdots(x-e_5^{})}
\right\} \Psi
\end{equation}
to a hypergeometric equation and, in fact, explicitly demonstrated
the first nontrivial example of integrability of equation
(\ref{fuchs}). Note that  even number of regular singular points is
essential in the context of uniformization by Whittaker's approach
\cite{whittaker1}, as hyperelliptic curves have always even number
of branch points.

\subsection{Where do hypergeometric equations come from?}
Whittaker does not elucidate the nature of his
conjecture\footnote{We were unsuccessful in finding any unpublished
manuscript material of Edmund Whittaker on his conjecture in
Edinburgh and London. Nevertheless the author  thanks the staff of
the London Mathematical Society Archives for their help in seeking
such information.} (\ref{conj}) but  the idea (not formulated) goes
back  to H.\,Weber (see formulae{} (10--16) in \cite{weber}),
although Weber considered a conformal representation of
multi-connected areas. The change of variables $x \mapsto y$ and the
substitution
\begin{equation}\label{subs}
\Psi(x)={\textstyle\sqrt{\mathstrut y\,}}\,\, \widetilde\Psi(y)\,,
\end{equation}
lead to the appearance of the  hypergeometric equation \cite{weber, whittaker2}
\begin{equation}\label{hyper}
y\,(y-1)\,\widetilde\Psi_{{\mathit{yy}}}+ \big[(a+b+1)\,y
-c\big]\,\widetilde\Psi_{y}+a\,b\,\widetilde\Psi=0\,.
\end{equation}
Such  a reduction is not a common property of hyperelliptic curves
and  the relation between equations (\ref{fuchs}) and (\ref{hyper})
(when it exists) depends on the kind of substitution. Nevertheless,
the question arises concerning Burnside's formulae (\ref{burnside}).
An attempt to apply such an argument (including the substitution
(\ref{subs}) to the curve (\ref{bs}) was undertaken in \cite{shabde}
but the conjecture certainly does not fit our example, because the
group contains a parabolic element: $x(\tau+4) = x(\tau)$. We
observe however that the $Q$-functions (\ref{conj}) and (\ref{Ds2})
differ from each other by only the numeric multiplier $\frac 43$.
Without taking into account this multiplier, all the accessory
parameters $A_{1,2,3}$ in Fuchsian equations (\ref{fuchs}) are equal
to zero (the same zero as in formula (\ref{Ds2})). The examples of
Burnside (1) and Whittaker (\ref{whit}) are thus the simplest ones
with known monodromies.

An explicit connection of the example in question with the modular
group and the hypergeometric equation  is given by the formula
(\ref{J}).  Indeed,  function $J(\tau)$ satisfies the Fuchsian
equation of hypergeometric type (\ref{hyper}) with
$\big(a=b=\frac{1}{12},\;c=\frac23\big)$ and therefore $J$ is
defined as inversion of a quotient of its two solutions. We do not
display here numerous forms of such representations (see for example
Klein's formulae \cite[p.\,61]{klein} or formulae (22)--(25) in
\cite[14.6.2]{erdei}). From (\ref{J}) we deduce that $x^4=z$ is a
root of the polynomial
$$
(z^2+14\,z+1)^3- 108\,z\,(z-1)^4\!\cdot\! J(\tau)=0\,.
$$
Hence by this construction $\sqrt[\leftroot{1}\uproot{1}4]{z}(\tau)$
is a {\em globally single-valued\/} function no matter which branch
of this root is chosen. This is a hidden form of Burnside's function
$x(\tau)=\sqrt[\leftroot{1}\uproot{1}4]{z}(\tau)$ in
(\ref{burnside}).

There is a simpler relation and a reduction of the Fuchsian equation
(\ref{last}) to the hypergeometric equation. Namely, the function
$\psi(z)=\Psi(\sqrt[\leftroot{1}\uproot{1}4]{z\,})$ satisfies, as it
is readily checked, the Fuchsian equation in $z$ with three regular
singularities whereupon we  can obtain another form of the integral
and  its series representations.

Analogous situation takes place for Whittaker's example
(\ref{whit})--(\ref{fuchs}). A simple reduction to the
hypergeometric equation immediately follows from the Fuchsian
equation for the second variable $y$.  Indeed, after the change of
variables
\begin{equation}\label{sub}
\big(x,\,\Psi(x)\big) \mapsto \big(y,\,
\widetilde\Psi(y)\big):\qquad \widetilde\Psi(y)={\textstyle
\sqrt{y_x^{}}}\,\Psi(x)\,,
\end{equation}
(this is not  the Weber--Whittaker change (\ref{subs})) to the
canonical form of second order equation, we readily get an equation
\begin{equation}\label{here}
\widetilde\Psi_{\mathit{yy}}=-\frac{6}{25}\,\frac{y^2+3}{(y^2-1)^2}\,
\widetilde\Psi\,,
\end{equation}
and its solution, say, in terms of hypergeometric series:
$$
\begin{array}{l}
\ds
\widetilde\Psi_1(y)=\qquad\quad\;\;(y^2-1)^{\frac25}_{\mathstrut}
\cdot
{}_2^{}F_1^{}\!\!\left(\frac25,\,\frac15;\,\frac45\,\Big|\,\frac{1-y}{2}
\right)_{{}_{\ds\mathstrut}}^{},\\
\ds\widetilde\Psi_2(y)=(y-1)^{\frac15}_{\mathstrut}
(y^2-1)^{\frac25}_{\mathstrut}\cdot
{}_2^{}F_1^{}\!\!\left(\frac35,\,\frac25;\,\frac65\,\Big|\,\frac{1-y}{2}
\right)^{\ds\mathstrut}.
\end{array}
$$
The monodromy group of the hypergeometric equation (\ref{hyper}) is
always Fuchsian if the difference of exponents is the reciprocal of
an integer 1, 2, \ldots, etc. In the case (\ref{here}) the Fuchsian
exponents for all the points $y=\{1,\;-1,\,\infty \}$ are
$$
\Big\{\frac25,\,\frac35\Big\}\,,
$$
so that  the monodromy groups for  Whittaker's equations
(\ref{fuchs}) and (\ref{here}) have both genus zero. Their
intersection is a genus 2 subgroup, i.\,e. group with four
generators uniformizing  equation (\ref{whit}). In fact, this was
done by Whittaker \cite{whittaker2}.

\subsection{The conjecture and hypergeometric reducibility}
Whittaker's conjecture admits an additional treatment. Let us write
equation of the curve (\ref{hyp}) in the form  $y^2=x^{2g+1}+E(x)$,
where $E(x)$ is a polynomial of degree $2\,g$. Bring Whittaker's
$Q(x)$-function (\ref{fuchs}) into the form
\begin{equation}\label{whit2}
Q(x)=-\frac38\left\{\frac{f_x^2}{f^2}-
4\,g\,(g+1)\,\frac{x^{2g-1}}{f}+\frac{E''+A(x)}{f} \right\}\;\,
\end{equation}
and compare (\ref{whit2}) with his conjecture (\ref{conj})
\begin{equation}\label{conj2}
Q(x)=-\frac38\left\{\frac{f_x^2}{f^2}-
4\,g\,(g+1)\,\frac{x^{2g-1}}{f}+\frac{2\,g+2}{2\,g+1} \frac{E''}{f}
\right\}\,.
\end{equation}
It immediately follows that accessory polynomial is a symmetrical
polynomial in branch points $e_j^{}$, i.\,e. $(2\,g+1)\,A(x)=E''$.
(The polynomial $A(x)$, as symmetrical one in $e_j^{}$, was
considered in work \cite{whittaker3}). When is it possible to reduce
equation (\ref{fuchs}) to the hypergeometric one? One way is to make
use of the substitution (\ref{sub}), whereupon we get
$$
\widetilde\Psi_{\mathit{yy}}=-\frac14\, \frac{4g(g+1)(2\,g+1)^2y^6-
(2g+1)\big(8\,x^3E'''+3(2g+1)\,x^2\,A(x)+
\cdots\big)\,y^4+\cdots}{\big((2g+1)\,y^2+x{\ds E'}-(2g+1)\,E
\big)^4}\, \widetilde\Psi\,.
$$
If we take into account that this rational function has not to
depend on $x$, we obtain from the denominator that $E(x)=a$, where
$a$ is a constant. One can give a rigorous form to this reasoning
but we speak here about only motivation. We  get
$$
\widetilde\Psi_{\mathit{yy}}=-\frac14\,
\frac{4\,g\,(g+1)\,(y^2+3\,a) - 3\,x^2\,A(x)}
{(2\,g+1)^2\,(y^2-a)^2}\,\widetilde\Psi\,,
$$
so that the only possibility is to put the accessory polynomial
$A(x)$ equal to zero
\begin{equation}\label{back}
\widetilde\Psi_{\mathit{yy}}=-
\frac{g\,(g+1)}{(2\,g+1)^2}\,\frac{y^2 + 3\,a}{(y^2-a)^2}\,
\widetilde\Psi
\end{equation}
and the monodromy automatically becomes   Fuchsian. Moreover it
becomes triangle group with angles $\frac{2\,\pi}{2g+1}$. This is
nothing else but all the examples of the papers
\cite{whittaker2,whittaker3,dhar} and elementary treatment of the
conjecture from the standpoint of such a reduction.

We may reverse the reasoning. Let us make the hyperelliptic change
of variables $y\mapsto x$ by the rule $y^2=x^{2g+1}+E$ and the
reverse substitution (\ref{sub}) into the hypergeometric equation
(\ref{back}) with the monodromy {\em known to be  Fuchsian\/}. We
shall obtain new Fuchsian equation in variable $x$. Require that
this equation has singularities at hyperelliptic branch points
$x=e_j^{}$ and coincides with Whittaker's form (\ref{whit2}). What
accessory polynomial $A(x)$ does satisfy these requirements? The
answer, as it follows from the previous constructions, is the
Whittaker conjecture (\ref{conj2}).

The conclusion is not changed if we consider the parabolic version
$$
Q(x)=-\frac12\left\{\frac{f_x^2}{f^2}-\frac{2\,g+2}{2\,g+1}\,
\frac{f_{\mathit{xx}}^{}}{f}\right\},
$$
but Burnside's example does not lead to a simple reduction of the
type $x \mapsto y$ (see (\ref{Qy})), although it admits  another
one. See \cite{br2} for details on the transformations of the type
(\ref{sub}) and additional discussion  the conjecture in
\cite{chud}. The conjecture is not true in general \cite{chud} but
the question ``when and why does it work?'' remains still open. Yet
another and perhaps simplest motivation to the hypergeometric
reducibility is the fact that equations (\ref{whit}) and
$y^2=x^{2g+1}+1$ represent an algebraic function
$x(y)=\!\!\!\sqrt[\leftroot{-4}\uproot{3}2g+1]{y^2-1}$ with three
branch points $\{-1,\,1,\,\infty \}$. But three points,
independently of kind of uniformization (Whittaker's or parabolic),
always lead to a hypergeometry ${}_2F_1$. Correlation between the
conjecture and equations with three branch-points was also remarked
in work \cite{chud}.

\section{Uniformization and covers  of  tori}

If a curve covers an elliptic torus then an integrable Fuchsian
equation on  torus is naturally to be expected.  This is the topic
of the present section.

Riemann surfaces of higher genera have negative curvature. Therefore
the simplest way to get nontrivial Fuchsian equations and
uniformization over tori is to consider torus with at least one
puncture. Such a problem is described by equation
$\Xi_{\alpha\alpha}=\big(\!-\!\frac14\wp(\alpha)+A \big)\,\Xi$ which
was already considered in the literature. First paper on this topic
was the work \cite{keen} (see also subsequent works of L.\,Keen)
however, up till now, no example with explicit analytic formulae has
been obtained.  In addition to the lemniscate  and equi-unharmonic
cases found in \cite{keen}, brothers Chudnovsky revealed
\cite{chud2} two more exceptional cases when an accessory parameter
is known and associated Fuchsian group has an arithmetic-algebraic
nature. See their work \cite{chud2} for more references and
additional discussion to this problem. Solutions of Fuchsian
equations, the problem of inversion and, what is important,
transformations between the equations are not considered in these
works. On the other hand, if we have correct accessory parameters
and solvable Fuchsian uniformization  for some algebraic curve and,
in turn, this curve can be realized as a cover over torus, then we
could automatically obtain non-trivial  examples on torus. The
general mode of getting such results is as follows.

Let we have a Fuchsian equation  in variable
$x\in\overline{\mathbb{C}}$, i.\,e. Fuchsian equation on plane, say,
(\ref{fuchs}). Let now $R(x,\alpha)=0$ be a formula for the cover
wherein $\alpha$ is the global parameter on torus being covered. In
most general case such formulae have the form
\begin{equation}\label{R}
\widetilde R\big(x,\wp(\alpha),\wp'(\alpha)\big)=0
\end{equation}
with polynomial function $\widetilde R$. Considering equation
(\ref{R}) as a transcendental (non algebraic) change of variables
$x\mapsto \alpha$ we can transform the initial Fuchsian
``$x$-equation'' into an equation in $\alpha$ (``$\alpha$-equation).
If ``$x$-equation'' was of Fuchsian class with correct accessory
parameters then the ``$\alpha$-equation'' will be of the same class.
Burnside's curve is just the case. In other words, from the
uniformization point of view, not merely punctured tori should be
searched for, but situations when the tori are covered by nontrivial
algebraic curves.

Since the examples that follows are the first exactly solvable ones
along these lines we expound all this at greater length in
Sec.\,IV.C.

\subsection{Cover of torus}
First (non-complete but hyper- and non-hyperelliptic) examples of
cover of torus $g=1$ by curves of  higher genera were obtained in
the first volume of Legendre's {\em Trait\'e\/} \cite{legendre}. We
consider a hyperelliptic example of Legendre generalized by Jacobi
\cite[{\em Werke} {\bf I}: pp.\,375--381]{jac}\footnote{In the same
place on p.\,377 Whittaker's curve (\ref{whit}) appeared. In this
section we return to the standard conventions for Weierstrassian
$\zeta,\wp$-functions: $\wp(\alpha)=\wp(\alpha|\om,\omp)$,
$e=\wp(\om|\om,\omp)$, etc.}:
\begin{equation}\label{AB}
y^2=x\,(x-1)(x-a)(x-b)(x-a\,b)\,.
\end{equation}
Jacobi found a substitution (simpler version of Legendre's one
\cite[{\bf I}: p.\,259]{legendre}) of second degree $x \mapsto
\lambda$:
\begin{equation}\label{jacobi}
\lambda=\frac{(1-a)(1-b)\,x}{(x-a)(x-b)}\,,
\end{equation}
which leads to the fact that both  holomorphic differentials for the
curve (\ref{AB}) reduce to elliptic differentials for the tori
\begin{equation}\label{k}
(\lambda,\,\mu)^{}_{\pm}:\qquad
\mu^2=\lambda\,(\lambda-1)\,(k\,\lambda-1)\,,\qquad
k=k_{\pm}^{}=-\frac{\big(\sqrt{\mathstrut a\,}\pm\sqrt{b\,}\big)^2}
{(a-1)(b-1)}\,.
\end{equation}
We transform this equation into the canonical Weierstrassian form
$$
w^2=4 \Big(z-\mfrac{2\,k-1}{3}\Big)\!
\Big(z-\mfrac{2-k}{3}\Big)\!\Big(z+\mfrac{k+1}{3}\Big)
$$
with the help of obvious scale transformation and subsequent
parametrization:
$$
\Big(z=k\,\lambda-\frac{k+1}{3}\,,\quad w=2\,k\,\mu\Big) \quad
\Rightarrow\quad
\big(z,\,w\big)=\big(\wp(\alpha),\,\wp'(\alpha)\big)\,.
$$
The Burnside curve corresponds to the following parameters
$$
a=-1,\quad b=i\,,\quad k_{\pm}^{}=\frac{1\pm\sqrt{2\,}}{2}\,,\quad
g_2^{}=\frac53\,, \quad g_3^{}=\mp\,\frac{7}{27}\sqrt{2\,}
$$
and  the two tori are isomorphic to the one classical torus with a
complex multiplication:
$$
J\Big(\mfrac{\omp}{\omega}\Big)=\frac{5^3}{3^3}\,.
$$
It has another standard form  $\wp'^2=4\,\wp^3-30\,\wp-28$ due to
formula
$$
\wp\!\left(\!z;\mfrac53,\mfrac{7\sqrt{2}}{27}\right)=\frac{\sqrt{2}}{6}\,
\wp\!\left(\!\sqrt{\textstyle\!\frac{\sqrt{2}}{6}}\,z;30,28\right).
$$
The values of corresponding periods are computed with use of modular
inversion problem:
$$
\omega=\sqrt[\uproot{3}4]{\mfrac{12}{5}
\,g_2^{}\big(\!\sqrt{2}\,i\big)}\,,
\qquad\omp=\frac{i}{\sqrt{2}}\,\omega\,.
$$
Using some properties of modular functions one can show \cite{br3}
that this Burnside's torus has an exact solution to the constant
$\omega$ since $g_2^{}\big(\sqrt{2}\,i\big)=\frac53\,\pi^4\,
\widehat\eta^8\big(\sqrt{2}\,i\big)$. We get
$$
\begin{array}{c}
\,\omega=\pi\,\sqrt{2}\,\widehat\eta^2\big(\sqrt{2}\,i\big)=
2.118\,156\,723\,947\,863\,188\,505\,038\,347\,005\,72\!\ldots\,,\\\\
\omp=\!\!\sqrt[-2]{2}\,i\!\cdot\!\om\,,\qquad \ompp=-\om-\omp\,,
\end{array}
$$
where, to avoid confusion between Weierstrass's and Dedekind's
standard notations for their $\eta$-functions, we denoted Dedekind's
one as $\widehat\eta(\tau)$: $\widehat\eta(\tau)=e^{\frac{\pi
i}{12}\tau}_{\mathstrut} \prod\!\big(1-e^{2\pi i n\tau} \big)$.

Assuming all the introduced parameters to be fixed and taking upper
sign in $k_{\pm}$, we view the Jacobi substitution (\ref{jacobi}) as
an explicit 2-sheeted cover of the torus
\begin{equation}\label{tor}
\wp'(\alpha)^2=4 \Big( \wp(\alpha)-\mfrac{\sqrt{2\,}}{3}\Big) \!
\Big(\wp(\alpha)+\mfrac{3+\sqrt{2\,}}{6} \Big)\! \Big(
\wp(\alpha)-\mfrac{3-\sqrt{2\,}}{6} \Big)
\end{equation}
by the $x$-planes or a fundamental 10-gon for the function $x(\tau)$
in $\tau$-plane \cite{br2}.  More precisely, the substitution
(\ref{jacobi}), i.\,e. equation (\ref{R}), has the form
\begin{equation}\label{cov}
\wp(\alpha)={\textstyle e'}+
\frac{1+3\,e}{x(\tau)-i}-i\,\frac{1+3\,e}{x(\tau)+1}\,,
\end{equation}
where the quantities $(e,e',e'')$ are taken from (\ref{tor}) in that
order. The formula (\ref{cov}) constitutes another representation of
a Riemann surface defined by the Burnside curve or, which is the
same, an equivalent {\em transcendental\/} representation
$R(\alpha,\tau)=0$ of the curve (\ref{bs}) itself in terms of
meromorphic functions on covering and covered  surfaces. Omitting
argument $\tau$ in  (\ref{cov}) we shall deal with {\em
algebraic-transcendental\/} version $R(\alpha,x)=0$ of this
representation.

\subsection{Structure of the cover (\ref{cov})}
Let us consider   ramification schemes of  equation (\ref{cov}). The
ramification points $\alpha_{j}^{}$ of this representation, as a map
$\alpha \mapsto x$, are determined from the equation
$$
\big(6\,\wp(\alpha)-3+\sqrt{2\,}\big)\!
\big(6\,\wp(\alpha)+21+13\,\sqrt{2\,}\big)=0
$$
(this is discriminant of the formula (\ref{cov})). Therefore its
solutions have the form
$$
\begin{array}{l}
\big(\alpha_{1}^{}=\ompp, \;\; x=- i\,\sqrt{i\,}\big), \qquad
\big(\alpha_{2,3}^{}=\pm\, \Kappa, \;\; x=i\,\sqrt{i\,}\big),
\end{array}
$$
where $\pm \Kappa$ are solutions of the
transcendental equation
$$
\wp(\Kappa)=-\frac72-\frac{13}{6}\,\sqrt{2\,}\qquad
\left(\Rightarrow
\wp'(\Kappa)=\sqrt[\leftroot{1}4]{32}\,\big(7+5\sqrt{2}\big)\,i
\right).
$$
At the point $\alpha_1^{}=\ompp$ we have two independent holomorphic
series:
$$
x(\alpha)=\frac{1-i}{2}\,\Big(\sqrt{2\,}  \pm \sqrt[4]{2\,}\,
(\alpha-\ompp)+\frac12\,(\alpha-\ompp)^2 \pm \cdots\Big)
$$
and hence the ramification scheme is \{1,\,1\} that means two
non-permutable holomorphic branches. At  each point $\pm \Kappa$ the
ramification scheme is \{2\} (two permutable holomorphic branches).
This follows obviously from the expansions
$$
x(\alpha)=i\,\sqrt{i\,}+\sqrt[\leftroot{1}4]{32\,}\,\sqrt{\alpha-\Kappa\,}+
\cdots, \qquad
x(\alpha)=i\,\sqrt{i\,}-i\,\sqrt[\leftroot{1}4]{32\,}\,\sqrt{\alpha+\Kappa\,}+\cdots.
$$
Ramification scheme at infinity $\wp(\alpha_4^{})=\infty$ is trivial
since there is no  ramifying here in fact:
$(x-i)(x+1)\sim\alpha^2+\cdots$. A check of the Riemann--Hurwitz
formula
\begin{equation}\label{hur}
\widetilde g=\frac12\,{\sum}_j \big(q_j^{}-1\big)+N\,(g-1)+1
\end{equation}
gives
$$
\widetilde g=\frac12\,\Big((1-1)+(2-1)\,2\Big)+2\,(1-1)+1=2
$$
as it must. Here $\widetilde g$ and $g=1$ are the genera of the
cover and surface under covering respectively, $N=2$ is number of
sheets of the cover (\ref{cov}), and $q_j^{}$ are indices of the
ramification at all the branch points $\alpha_{1,2,3}^{}$.

In the back direction $x\mapsto\alpha$ the representation
(\ref{cov}) is also  $(1 \mapsto 2)$-map since two copies of torus
cover the $x$-plane. Ramification schemes are as follows
$$
\begin{array}{rl}
\wp'(\alpha)=0:&\quad\left\{
\begin{array}{llc}
\big(\alpha_1^{}=\omega, & \!\!x=\{1,\,-i\}\big), &\quad\Big\{\{2\},
\{2\}\Big\}_{\ds\mathstrut}\\
\big(\alpha_2^{}=\omp, &\!\! x=\{0,\,\infty\}\big), &
\quad \Big\{\{2\}, \{2\}\Big\}_{\ds\mathstrut}\\
\big(\alpha_3^{}=\ompp, &\!\! x=-i\,\sqrt{i\,}\big),&
\quad \{1,\,1\}^{\ds\mathstrut}\\
\end{array}
\right.,\\ \\
\wp'(\alpha)=\infty:&\qquad\; \big(\alpha_4^{}=0,\,\;\, x=\{-1,\,i
\}\big), \quad\;\; \Big\{\{2\}, \{2\}\Big\}\,.
\end{array}
$$
For example in the neighborhood of the first branch point
$\big(\alpha=\omega,\;x=1\big)$ we compute formally Puiseux series
\begin{equation}\label{puiseux}
\alpha(x)=\omega-i\,\textstyle{\sqrt{i\,\sqrt{2\,}+i\,}}\,
\sqrt{x-1}+
\mfrac{i\,\sqrt{26\,\sqrt{2\,}+34}-\sqrt{26\,\sqrt{2\,}-14}}{24}\,
\sqrt{x-1}^{\,3}+\cdots.
\end{equation}
These series are important and discussed in Sec.\,V.D. Checking the
formula (\ref{hur}) under $g=0$ and $N=2$ (order of the function
$\wp(\alpha)$) we obtain again
$$
\widetilde g=\frac12\,\Big((2-1)\,4 +(1-1)\,1 +(2-1)\,2
\Big)+2\,(0-1)+1 =2\,.
$$

\subsection{Solvable Fuchsian equations on torus}
A direct consequence of the preceding arguments is a linear ordinary
differential equation  of Fuchsian type on the torus (\ref{tor})
with well defined accessory parameters. Indeed, making use of
Jacobi's substitution (\ref{jacobi}) and changing $\Psi$ in
(\ref{last}) by the rule
$$
\Psi\mapsto\Lambda:\qquad
\Lambda(\lambda)=\sqrt{\lambda_x}\,\Psi(x)=
\frac{\sqrt{x^2+i\,}}{(x+1)(x-i)}\,\Psi(x)\,,
$$
we arrive at equation with five regular singularities
$\lambda_j^{}=\big\{0,\,1,\,-2\pm2\,\sqrt{2\,},\,\infty \big\}$:
\begin{equation}\label{Lambda}
\begin{array}{l}
\Lambda_{\lambda\lambda}^{}=\ds
\;\;\,\frac12_{\ds\mathstrut}^{}\,Q(\lambda)\, \Lambda={}\\
\phantom{\Lambda_{\lambda\lambda}^{}}=\ds -\frac14\,
\frac{\lambda^6+4\,\lambda^5+16\,\lambda^4-56\,\lambda^3+68\,\lambda^2
-48\,\lambda+16}{\lambda^2\,(\lambda-1)^2\,(\lambda^2+4\,\lambda-4)^2}^
{\ds\mathstrut}_{}\, \Lambda\,.
\end{array}
\end{equation}
This equation is remarkable in itself because not every algebraic
change of variables (formula (\ref{jacobi}) in our case) in Fuchsian
equation with rational coefficients and Fuchsian monodromy leads to
equation, again, with rational coefficients and Fuchsian monodromy.
Singularities $\lambda_j^{}=-2\pm2\,\sqrt{2\,}$ correspond to
elliptic edges of the second order and
$\lambda_j^{}=\{0,\,1,\,\infty\}$ are parabolic singularities. See
\cite{br2} for further application of the Fuchsian equation
(\ref{Lambda}) to uniformization.

Considering the parametrization of the torus (\ref{k})
$$
\lambda=2\,\big(\sqrt{2\,}-1\big)\,\wp(\alpha)+\frac{2\,\sqrt{2\,}-1}{3}
$$
as a subsequent change $\lambda\mapsto\alpha$ and supplementing with
the formula
$$
\Lambda\mapsto\Xi:\qquad
\Xi(\alpha)=\textstyle{\sqrt{\alpha_{\lambda}^{}{\phantom{\mathstrut}}}}\,\Lambda(\lambda)
=\mathrm{const}\cdot\!\!
\sqrt[\leftroot{-2}\uproot{2}-2]{\wp'(\alpha)}\,\Lambda(\lambda)\,,
$$
we get the   equation
\begin{equation}\label{Q}
\Xi_{\alpha\alpha}=\frac12\,Q(\alpha)\,\Xi\,,
\end{equation}
where function $Q(\alpha)$ (it must be elliptic one) is given by the
following expression:
$$
Q(\alpha)=-\big\{\lambda,\,\alpha
\big\}+\lambda_\alpha^2\,Q(\lambda)=
6\,\wp(2\,\alpha)+4\,\big(\sqrt{2\,}-1\big)^2\wp'(\alpha)^2\,
Q(\lambda)\,.
$$
Carrying out some simplifications we obtain the sought-for result.

{\bf Proposition\;3.}  {\em Fuchsian  equation on torus
corresponding to Burnside's paramet\-ri\-za\-tion
{\rm(\ref{burnside})} is a linear ordinary differential equation on
torus {\rm(\ref{tor})} with five regular singularities. The equation
and accessory parameters have the form\/}
\begin{equation}\label{new}
\begin{array}{l} \ds\Xi_{\alpha\alpha}=\ds
\Big\{\!\!-\!
\mfrac{1^{\mathstrut}}{4}\big(\wp(\alpha)+\wp(\alpha-\omega)+
\wp(\alpha-\omp)\big) -\mfrac{3}{16}
\big(\wp(\alpha-\Kappa)+\wp(\alpha+\Kappa)\big)+{}
\\\\
\ds \phantom{\ds\Xi_{\alpha\alpha}=
\Big\{}{}\!+\mfrac{9}{64}\sqrt[\leftroot{1}4]{8\,}\,i\,\big(
\zeta(\alpha-\Kappa)-\zeta(\alpha+\Kappa)\big)+ \mfrac{9}{32}
\big(\sqrt[\leftroot{1}4]{8\,}\,i\,\zeta(\Kappa)+2\sqrt{2}+2\big)\Big\}
\,\Xi\,.
\end{array}
\end{equation}

{\em Remark\;$2$\/}. Accessory parameters in Fuchsian ODE  on torus
(\ref{Q}) are by definition coefficients in front of $\alpha^{-1}$
in function $Q(\alpha)$. These are multipliers of $\zeta$-functions
in (\ref{new}) and a free term.

In so far as Burnside's example is a  curve with maximal symmetries,
this equation provides, perhaps, the simplest example of
non-triviality: solvable Fuchsian {\sc ode} on  torus with  Fuchsian
monodromy, solution, explicitly known inversion (see Sec.\,V.D), and
algebraic curve. Solution of equation (\ref{new}), in many forms, is
given by back transformations written above. Invoking solutions of
equation (\ref{last}) described in \cite{br2} we readily obtain one
of such forms:
$$
\Xi(\alpha)=\sqrt{\alpha_x}\,\Psi(x)={\textstyle
\sqrt{\big(x-i\,\sqrt{i\,} \big)\,y}}\cdot\big(A\,K(x^2)+B\,K'(x^2)
\big)\,,
$$
where $K$ and $K'$ are complete elliptic integrals \cite{erdei} and
the quantities $x,x^2,y$, as functions of $\alpha$, have to be
expressed from the following pair of algebraic equations
(consequence of the equations (1), (\ref{tor}), and (\ref{cov})):
\begin{equation}\label{Palpha}
x^2=(i-1)\,\frac{\wp(\alpha)+e'-2\,e} {\wp(\alpha)-e'}\,x+i\,,\qquad
\wp'(\alpha)=-\frac{6\,e+2}{\sqrt{1+i\,}}\, \frac{\big(x\pm
i\,\sqrt{i\,}\big)\,y}{(x-i)^2(x+1)^2}\,.
\end{equation}
The ultimate formula can be represented in various simplified forms
but solutions of all the Fuchsian equations we consider are
essentially multi-valued functions. This indispensable property
reflects the required character of their monodromy groups.

{\bf Corollary.} {\em The quotient
\begin{equation}\label{Xi}
\tau=\frac{\Xi_2^{}(\alpha)}{\Xi_1^{}(\alpha)} \qquad
\Leftrightarrow\qquad \mbox{\rm (\ref{cov})\quad and\quad
(\ref{last})}
\end{equation}
is the global parameter on both the Burnside Riemann surface  and
corresponding orbifold $\boldsymbol{\mathfrak{T}}$, defined by
monodromy group of the equation {\rm (\ref{new})}. Inversion of the
ratio {\rm(\ref{Xi})} is, by construction, a globally single-valued
analytic function $\alpha=\alpha(\tau)$\/}.

Another form of  equation (\ref{new}) can be useful:
$$
\!\begin{array}{l} \ds\Xi_{\alpha\alpha}=\ds
\Big\{\!-\mfrac{1}{4}\big(\wp(\alpha)+\wp(\alpha-\omega)+
\wp(\alpha-\omp)\big)-{}\\\\
\ds \phantom{\ds\Xi_{\alpha\alpha}= \Big\{\!} -\mfrac{3}{16}
\Big(\wp(\alpha-\Kappa)+\wp(\alpha+\Kappa) +
\mfrac{3\,(7+5\sqrt{2})}{\wp(\alpha)-\wp(\Kappa)} -3\,
\sqrt{2}-3\Big)\Big\}\,\Xi\,.
\end{array}
$$
By renormalizing the global parameter on torus
$\alpha\mapsto\tilde\alpha$ by the rule $\alpha=\omp\,\tilde\alpha$
and bringing Weierstrassian functions $\wp(z|\omega,\omega')$ into
the form $\wp(z|1,\mu)=:\widetilde\wp(z|\mu)$ with one modulus:
$$
\zeta(\alpha|\om,\omp)=\frac{1}{\omp}\,
\widetilde\zeta\big(\tilde\alpha|\sqrt{2}\,i\big),\qquad
\wp(\alpha|\om,\omp)=\frac{1}{\omp^2}\,
\widetilde\wp\big(\tilde\alpha|\sqrt{2}\,i\big)
$$
we can present  equation (\ref{new}) in another canonical form:
$$
\begin{array}{l} \ds\Xi_{\tilde\alpha\tilde\alpha}^{}=\ds
\Big\{\!- \mfrac{1^{\mathstrut}}{4}\big(\widetilde\wp(\tilde\alpha)
+ \widetilde\wp(\tilde\alpha-1)+
\widetilde\wp(\tilde\alpha-\sqrt{2}\,i)\big) -\mfrac{3}{16}
\big(\widetilde\wp(\tilde\alpha-\tilde\Kappa)
+\widetilde\wp(\tilde\alpha+\tilde\Kappa)\big)-{}
\\\\
\ds \phantom{\ds\Xi_{\alpha\alpha}= \Big\{}{}-\mfrac{9}{64}\,M\big(
\widetilde\zeta(\tilde\alpha-\tilde\Kappa)-
\widetilde\zeta(\tilde\alpha+\tilde\Kappa)\big)- \mfrac{9}{32}\,
M\big(\widetilde\zeta(\tilde\Kappa)+
(\!\sqrt[-2]{2}+1)\,M\big)\Big\} \,\Xi\,,
\end{array}
$$
where $M:=\sqrt[\leftroot{1}\uproot{1}4]{2\,}\,\omega=
\sqrt[\leftroot{1}\uproot{1}4]{8\,}\,\pi\,\widehat\eta^2(\sqrt{2}\,i)$.
Note that all the accessory parameters are the real quantities as
$\widetilde\zeta(\tilde\Kappa)=3.83102282421\!\ldots$ and this
equation, under real $\alpha$, is defined by a real even function of
$\alpha$.

Clearly, the points $\alpha_j^{}=\big\{0,\,\omega,\,\omp \big\}$
correspond to three punctures on the torus, because in the vicinity
of these points we have  expansions of the form
$$
\frac12\,Q(\alpha)=-\frac14\,\frac{1}{(\alpha-\alpha_j^{})^2}+\,\cdots\,.
$$
The  remaining points $\pm \Kappa$ correspond (locally) to elliptic
edges of  second order:
$$
\frac12\,Q(\alpha)=-\frac{3}{16}\,\frac{1}{(\alpha \mp \Kappa)^2}
\pm\frac{9}{64}\,\sqrt[\leftroot{1}4]{8\,}\,i\,\frac{1}{(\alpha \mp
\Kappa)}+\cdots\,.
$$
The {\sc Fig.\,1} illustrates these remarks. First example of
Poincar\'e's metric on the ``toroidal'' orbifold
$\boldsymbol{\mathfrak{T}}$ is presented in the work \cite{br2}.

\begin{figure}[htbp]
\centering
\definecolor{gray}{gray}{0.85}
\centerline{\hspace{-3.5cm}\unitlength=0.8mm
\begin{picture}(40,32)
\put(0.34,5.4){{\color{gray} \rule{8.27em}{6.05em}}}
\thicklines
\put(0.6,4.9){\line(1,0){16.0}}\put(0.6,5.1){\line(1,0){16.0}}
\put(17.8,4.9){\line(1,0){16.3}}\put(17.8,5.1){\line(1,0){16.3}}
\thinlines
\put(28,5){\line(1,0){9}}
\multiput(-0.4,29.8)(1,0){34}{\scriptsize.}
\thicklines
\put(-0.1,5.6){\line(0,1){11.8}}\put(0.1,5.6){\line(0,1){11.8}}
\put(-0.1,18.6){\line(0,1){11.5}}\put(0.1,18.6){\line(0,1){11.5}}
\thinlines
\put(0,25){\line(0,1){8}} \multiput(33.5,4.8)(0,1){26}{\scriptsize.}
\put(0,18){\circle{1.2}} \put(0,5){\circle{1.2}}
\put(17.2,5){\circle{1.2}}
\put(32,1){\footnotesize$2\omega$}
\put(16.1,1){\footnotesize$\omega$}
\put(-7.5,29){\footnotesize$2\omp$} \put(-5,17){\footnotesize$\omp$}
\put(-3,1){\footnotesize $0$} \put(22,24){{\footnotesize
$(\alpha)$}}
\put(-1.33,21.4){$\scriptstyle\boldsymbol{\times}$}\put(2.5,21.5){\footnotesize
--\Kappa} \put(-1.33,7.3){$\scriptstyle\boldsymbol{\times}$}
\put(2.5,7.3){\footnotesize \Kappa} \put(58,26){\circle{1.2}}
\put(62,25){-- \ punctures}
\put(54,15){\footnotesize$\; -\Kappa = 2.604826300529...\,i$}
\put(54,5){\footnotesize$\;\phantom{-}\Kappa =
0.390699665709...\,i$}
\end{picture}
} \caption{Orbifold\, $\boldsymbol{\mathfrak{T}}$ with fundamental
group $\pi_1^{}(\boldsymbol{\mathfrak{T}})$  of rank 6
(Proposition\;7) defined by monodromy group of equation (\ref{new}).
Crosses ``$\boldsymbol{\times}$''  stand for elliptic edges of 2nd
order.}
\end{figure}
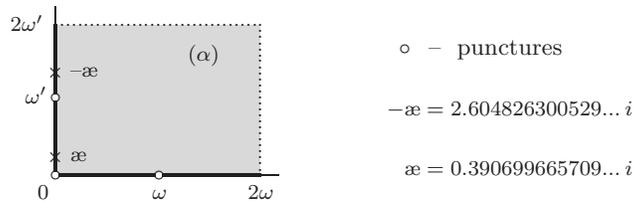

In order to obtain Fuchsian equation on torus, it is not necessary
for the torus to  be punctured. One may add  elliptic singularities
of finite order. Rankin \cite{rankin2} showed that Whittaker's
conjecture holds for Burnside's curve (1) and accessory parameters
$A_{1,2,3}$ in equation (\ref{fuchs})  have also zero values.
Therefore we can derive a Fuchsian equation on torus (\ref{tor})
corresponding to this kind of uniformization. Carrying out analogous
calculations, with the only difference that the previous
substitutions are applied to equation (\ref{fuchs}) with
$A_{1,2,3}=0$, we get yet another example of solvable Fuchsian
equation on torus.

{\bf Proposition\;4.} {\em The orbifold\/
$\boldsymbol{\mathfrak{T}}$ corresponding to Whittaker's
uniformization of Burnside's Riemann surface {\rm(\ref{cov})} is
defined by a Fuchsian equation on torus {\rm(\ref{tor})} with two
singularities and three accessory parameters. The equation has the
form}
\begin{equation}\label{new2}
\begin{array}{l} \ds\Xi_{\alpha\alpha}=
-\mfrac{3}{16}\Big\{ \wp(\alpha-\Kappa)+\wp(\alpha+\Kappa)
-\!\!\sqrt[\mbox{\tiny --}4]{32\,}\,i\!\cdot\!\zeta(\alpha-\Kappa)
+\!\!\sqrt[\mbox{\tiny
--}4]{32\,}\,i\!\cdot\!\zeta(\alpha+\Kappa)-{}
\\\\
\ds\phantom{\Xi_{\alpha\alpha}=\{--\;} -\big( \!\sqrt[\mbox{\tiny
--}4]{2\,}\,i\,\zeta(\Kappa)+\sqrt{2\,}+1\big)\Big\}\,\Xi
\end{array}
\end{equation}
{\em or, equivalently\/},
$$
\Xi_{\alpha\alpha}= -\mfrac{3}{16}\Big\{
\wp(\alpha-\Kappa)+\wp(\alpha+\Kappa)+
\mfrac{7+5\sqrt{2}}{\wp(\alpha)-\wp(\Kappa)}
-\sqrt{2\,}-1\Big\}\,\Xi\,.
$$

Explicit description of  monodromies (their genera, generators,
arithmetical properties, if any, etc) of $\boldsymbol{\mathfrak{T}}$
defined by equations (\ref{new}) and (\ref{new2}) is open question.
Erasing punctures in Fig.\,1 we obtain illustration of this
orbifold. Compared to equation (\ref{new}), the equation
(\ref{new2}) has only two singularities however the explicit
expression for inversion of the ratio (\ref{Xi}) for this equation
is unknown. This important point manifests itself in the fact that
{\em none of the uniformizing functions for any non-modular Fuchsian
equation  is known hitherto\/}.

\section{Analytical series}
 Power series are the widely exploited tools in the theory
of elliptic functions and have enormous number of applications
\cite{gannon}. In this section we develop corresponding technique
for the curve (1).

In order to obtain the series expansion of the function $x(\tau)$ we
may not make use of well-known Laurent's expansions of
Weierstrassian functions in so far as such expansions turn formally
$\wp$-functions in (\ref{burnside}) into modular forms of
fictitiously infinite weight. This is of course not the case.
Explanation is that  the pole $\tau=2$ of the function $x(\tau)$
lies on a singularity line (real axis) of $\wp$-functions in
(\ref{burnside}). The same complications occur for other branch
points of the curve (\ref{bs}): that is, roots of the polynomial
$x^5-x$. For  the same reasons the function (\ref{x45}) has no
seemingly singularities at points $x^4=5$. They are compensated by a
singularity of the function $\wp\big(\tau(x)|2,2\tau(x)\big)$.

To overcome this obstacle we should use Schwarz's equation
(\ref{Ds2}) satisfied by the function $x(\tau)$. The series
representation depend on  the location of  origin of the expansion.
It can lie inside the fundamental circle or on its border.

We shall use the natural notation $\om_j^{}$ for  pre-images of
branch points $e_j^{}$:
$$
x(\om_0^{})=\infty\,,\qquad x(\om_j^{})=e_j^{}\,, \qquad
e_j^{}=\big\{ 0,\,1, -1,\,i,-i\big\}\,.
$$

\subsection{Meromorphic derivative}
Let us define the {\em meromorphic derivative\/}
${\boldsymbol{\mathfrak D}}$ of  nonconstant function $x=x(\tau)$ by
the formula
$$
{\boldsymbol{\mathfrak D}}:\quad\big[ x,\,\tau\big] := \frac{\big\{
x,\,\tau\big\}}{x_\tau^2} \quad\Leftrightarrow\quad
\big[x,\,\tau\big] = -\big\{\tau,\,x \big\}\,.
$$
Motivation for introducing this object is the fact that
${\boldsymbol{\mathfrak D}}$ is a simplest and lowest order
differential combination of a meromorphic function on a curve which
is  meromorphic function as well. It also determines Schwarz's
equation (\ref{Ds1}). Properties of the object
${\boldsymbol{\mathfrak D}}$ follow from properties of the
Schwarzian:
\begin{equation}\label{prop}
\Big[X\big(q(\tau)\big), \, \tau\Big] = \left\{\!\big[X(q), \,
q\big]+\mfrac{1}{X_q^2}\, \big[q(\tau), \,
\tau\big]\!\right\}_{\!q=q(\tau)}\,.
\end{equation}

Since we shall deal with meromorphic functions on Riemann surfaces
(meromorphic automorphic functions, holomorphic and meromorphic
Abelian integrals), without loss of generality, by virtue of
(\ref{prop}), we may consider expansions only of the form
\begin{equation}\label{pole}
X=q^{-n}\,\big(A+ B\,q+C\,q^2+\cdots\big), \qquad n\in
\mathbb{Z},\;A\ne0\,.
\end{equation}
It follows that
\begin{equation}\label{poles}
\begin{array}{l}\ds
\big[ X,\, q \big ]=\frac{1}{A^2}\left\{\frac{1-n^2}{2\,n^2}+
\frac{n^2-1}{n^2}\,\frac{B}{A}\,q\right.+{}\\ \\
\ds \phantom{\big[ X,\, q \big ]=\frac{1}{A^2}}\;\left.{}+
\frac{2\,n\,(n^3-n-6)\,A\,C-3\,(n^4-n^2-2\,n+2)\,B^2}{2\,n^4\,A^2}\,q^2
+\cdots\right\}q^{2n}\,.
\end{array}
\end{equation}
This property entails the fact that ${\boldsymbol{\mathfrak D}}$ has
a pole/zero if, and only if the function $X(q)$ has a {\em fold\/}
zero/pole, i.\,e. $n\ne\pm1$, or a fold $a$-point (this is a point
$\tau$ where the function $x(\tau)$ takes the value $a$). Poles and
zeroes of ${\boldsymbol{\mathfrak D}}$ are always fold. The fold
$a$-points of analytic functions lead to violation of conformality
of analytic maps and, in the uniformization theory, correspond to
branch points of covers and singularities in Fuchsian equations.
Conversely, let $x(\tau)$ be meromorphic automorphic function on
some algebraic curve. Being analytic function of $\tau$, the
function $x(\tau)$ has only finitely many $a$-points. From
(\ref{prop})--(\ref{poles}) it follows that its
${\boldsymbol{\mathfrak D}}$-derivative is meromorphic function:
that is, rational function of $(x,y)$. This is nothing but equation
(\ref{Ds1}) without recourse to auxiliary Fuchsian linear {\sc ode}.

\subsection{Local and global parameters}
To determine the behavior of a local parameter $q=q(\tau)$ and,
therefore, the type of point, we set $x(\tau)=X(q)$. Let us apply
the object ${\boldsymbol{\mathfrak D}}$ to our main example. Let
$\om_j^{}$ be a zero or pole of the function $x-e_j^{}$. For example
the pole (\ref{pole}). Using the properties
(\ref{prop})--(\ref{poles}), an expansion of both sides of
(\ref{Ds2}) produces
$$
\begin{array}{r}
\ds \frac{1}{A^2}\left\{\Big(\frac{1-n^2}{2\,n^2} + \frac
BA\,\frac{n^2-1}{n^2}\,q+\cdots\Big)+ \frac{1}{n^2}\,\big[q,\,
\tau\big]\,q^2 \Big( 1-2\,\frac BA\,\frac{n-1}{n}\,q +
\cdots\Big)\right\}
q^{2n}\phantom{\,.}={}\\ \\
\ds =\frac{1}{A^2}\left\{-\mfrac{1}{2}+\mfrac BA\,q+\cdots\right\}
q^{2n}\,.\quad\;
\end{array}
$$
Balancing degrees with respect to $q$ entails
\begin{equation}\label{parabolic}
\big[q,\, \tau\big]\,q^2=-\frac12+\cdots.
\end{equation}
Integrating  this equation by series $q=a+b\,\tau+\cdots$ shows that
$a,\,b\ne 0$. Moreover, the inverse function $\tau(q)$ is never
meromorphic in $q$ as the ansatz $\tau=q^{-k}+\cdots$ involves
incompatibility with (\ref{parabolic}).  Hence, locally, $q$ has an
exponential behavior
\begin{equation}\label{dots}
q(\tau)=\exp\!\left(\frac{a\,\tau+b}{c\,\tau+d}\right)+\cdots \;
\mbox{holomorphic part}\,.
\end{equation}
Choice of local parameters has no restrictions except for condition
of being locally single-valued. Therefore we  may omit all the dots
in (\ref{dots}) and obtain the following representation:
\begin{equation}\label{abcd}
q(\tau)=\exp\!\left(\frac{a\,\tau+b}{c\,\tau+d}\right)\,,
\end{equation}
where $\big(a,\,b,\,c,\,d\big)$ are arbitrary constants.  Therefore
$\om_0^{}$ and also all the pre-images $\om_j^{}$ of other branch
points are  edges of parabolic cycles. We shall call the formula
(\ref{abcd}) {\em $q$-representation of the global coordinate
$\tau$\/}. Note that (\ref{abcd}) is the usual form for  local
parameters in the vicinity of such edges \cite[\S\,41]{ford} and the
order $n$ is not determined.

In the case of Whittaker's uniformization (\ref{fuchs}) we would
have $n=\pm 2$ and  an equation for the local parameter in the form
\begin{equation}\label{abcd'}
\big[ q,\, \tau \big ]=0+\cdots
\quad \Rightarrow \quad q(\tau)=
\frac{a\,\tau+b}{c\,\tau+d}+\cdots,
\end{equation}
so it is not necessary to change the global parameter $\tau$ to the
local $q$.  This is not surprising, because the pre-images
$\om_j^{}$ lie inside the fundamental polygon
\cite{whittaker1,whittaker2} and, as  is well-known, they are
Weierstrass's points, and the meromorphic function $x-e_j^{}$ on the
curve has a 2-nd order  ($n=\pm2$) pole/zeroes at the points
$\om_j^{}$.

Schwarz's equation (\ref{Ds1}), after the substitution (\ref{abcd}),
acquires  the form
\begin{equation}\label{X}
\big[ X,\, q \big]-\frac{1}{2\,q^2}\,\frac{1}{X_q^2}=Q(X,Y)
\end{equation}
with the $Q$-function having the following structure
$$
Q(X,Y)=\frac{\mu}{\big(X-e_j^{}\big)^2}+\cdots.
$$
We must have a restriction on $\mu$ because we look for meromorphic
solutions (\ref{pole}) of  equation (\ref{X}), so that $n$ has to be
integer:
$$
n^2\,(2\,\mu+1)=1\,.
$$
The $Q$-function for any Fuchsian equation with Fuchsian monodromy
must satisfy this relation.  If $\mu = -\mfrac12$ we arrive at the
parabolicity condition (\ref{parabolic})--(\ref{abcd}), i.\,e.
$n=\infty$. Otherwise we arrive at (\ref{abcd'}) and determine the
order $n$. For a Fuchsian group without parabolic edge at the point
$e_j^{}$ we obtain: 1) $\mu=-\frac38$ $(n=\pm 2)$: hyperelliptic
curves; 2) $\mu=-\frac49$ $(n=\pm3)$: curves where $x$-function has
a 3-rd order pole/zero (for example trigonal curves); etc.

\subsection{Integer $\boldsymbol{q}$-series for genus two (Burnside's curve)}
The general solution of equation (\ref{X}) is determined up to the
``exponential-linear'' substitution
$$
q \mapsto \exp\!\left( \frac{a\ln q+b}{c\ln q+d}\right).
$$
In order to have a meromorphic function we have to obtain analytic
Laurent series for the function $X(q)$ in the neighborhood of $q=0$.
The constants $\big(a,\,b,\,c,\,d\big)$ in (\ref{abcd}) cannot be
determined from the differential equation (\ref{X}). However, well
choice of these constants  would allow one to construct analogues of
the celebrated modular integer $q$-series but, in our case, they
would correspond to a nontrivial group of genus two. Running ahead
note that it is not evident at all that these series should be
integer. If so, we shall call such series {\em canonical
representation\/} (see remark\;4 further below).

The uniformness of the functions $X=x(\tau)$ and $Y=y(\tau)$ entails
the following ansatz for the polar expansion (\ref{pole})
$$
X=c_{\scriptscriptstyle-2}\,q^{-2}+
c_{\scriptscriptstyle-1}\,q^{-1}+
c_{\scriptscriptstyle0}+c_{\scriptscriptstyle1}\,q+\cdots.
$$
Appropriate normalization of the constants $(a,b,c,d)$ for this pole
is as follows:
$$
q=e^{\frac{\pi i}{4}\frac{2\tau-5}{\tau-2}}_{\mathstrut},\qquad
\tau\to 2+i\,0\,.
$$
Substituting this ansatz into (\ref{Ds2}) we get an integer series
indeed:
\begin{equation}\label{Xseries}
X= \frac12\,q^{-2}\, \big(1+2\,q^8-q^{16}
-2\,q^{24}+3\,q^{32}+2\,q^{40} -4\,q^{48} -4\,q^{56}+\cdots\big)\,.
\end{equation}
The general solution $x(\tau)$ of  equation (\ref{Ds2}) is obtained
after the substitution (\ref{abcd}).  The $q$-series for the second
uniformizing function $y(\tau)$ is derived from  equations
(\ref{Qy}), (\ref{X}) or the identity $Y^2=X^5-X$. We find that $Y$
is also determined by an integer series:
\begin{equation}\label{Yseries}
Y=\frac{\sqrt{2}}{8}\,q^{-5}\,\big(1-3\,q^8- 3\,q^{16}
+14\,q^{24}+6\,q^{32}-33\,q^{40}- 20\,q^{48}+\cdots\big)\,.
\end{equation}

Being the single-valued functions of $\tau$, the series
(\ref{Xseries})--(\ref{Yseries}) are the sought for canonical
representations of  functions (\ref{burnside}) in the vicinity of
their pole $\om_0^{}=2$ with correctly chosen constants
$(a,\,b,\,c,\,d)$. Analyzing Burnside's formulae and applying the
same technique, one obtains canonical expansions for all other
branch points $e_j^{}$:
\begin{equation}\label{tmp}
x(0)=0\,,\quad x\Big( \mfrac12\Big)=1\,, \quad x(\infty)=-1\,, \quad
x(1)=i\,, \quad x(-1)=-i\,.
\end{equation}

{\em Remark\;$3$\/}. Note  $x(1)=i$ rather than $1$  as it might
seem from (\ref{xt12}). The point $\tau=1\not\in\mathbb{H}^+$ and we
must use  a limiting passage to get correct value $x(1)$. Similarly,
$x(-1)\ne x(1)$ despite the formal fact that $x(-\tau)=x(\tau)$.
Transformation $\tau\mapsto -\tau$  preserves shape of function
(\ref{xt12}) but sends $\mathbb{H}^+\mapsto \mathbb{H}^-$ and
therefore does not belong to group $\mathrm{PSL}_2(\mathbb{Z})$. We
see that this sole point forbids such an automorphism as it must.

{\bf Proposition\;5.} {\em The canonical representations for
Laurent's series for Burnside's  functions {\rm (\ref{burnside})} do
exist. The developments are determined by the formulae {\rm
(\ref{Xseries})--(\ref{Yseries})} in the vicinity of the pole
$\om_0^{}=2$ and coordinate
\begin{equation}\label{local}
q=e^{\frac{\pi i}{4}\frac{\tau-2}{2\,\tau}}_{\mathstrut},\qquad
\tau\to 0+i\,0
\end{equation}
corresponds to canonical representation in the vicinity of the zero
$\om_1^{}=0$$:$
\begin{equation}\label{moon}
\left\{\!
\begin{array}{l}
X=e^{\frac34\pi i}_{\mathstrut}\,\quad
2\,q^2\,\big(1+2\,q^8+5\,q^{16}+10\,q^{24}+ 18\,q^{32}+
32\,q^{40}+\cdots\big)_{\mathstrut}\\
Y=e^{-\frac{1^{\ds\mathstrut}}{8}\pi
i}_{\mathstrut}\,\sqrt{2\,}\,q\,\big(1+9\,q^8+42\,q^{16}+147\,q^{24}+
444\,q^{ 32}+1206\,q^{40}+\cdots\big)
\end{array}
\right..
\end{equation}
The developments at branch points $e_j^{}=\{\pm 1,\,\pm i \}$ are
determined, up to  multipliers, by the one canonical series\/$:$
\begin{equation}\label{XY}
\left\{\!
\begin{array}{l}
X=\big\{\!\pm\!1,\,\pm i \big\}\,\big(1 + 4\,q^2 + 8\,q^4 + 16\,q^6+
32\,q^8 + 56\,q^{10} + 96\,q^{12}+ \cdots\big)\\ \\
Y=4\,\big\{\sqrt{\pm1},\,\sqrt{\pm i\,}\big\}\,q\, \big (1 +
6\,q^2+24\,q^4 + 80\,q^6+231\,q^8+ 606\,q^{10} +\cdots\big)
\end{array},
\right.
\end{equation}
where coordinates $q$, according to {\rm (\ref{tmp})}, have the form
$$
\begin{array}{l}
\Big\{q=e^{\frac{\pi i}{4}\frac{3\tau-2}{2\tau -1}}_{\mathstrut},
\quad\tau\to\mfrac12+i\,0\Big\}_{\mathstrut}\,,\\\\
\Big\{q=e^{\frac{\pi i}{4}(\tau-2)}_{\mathstrut}, \;\;\,\tau\to
+i\,\infty\;\;\,\Big\}\,,
\end{array}
\quad
\begin{array}{ll}
\Big\{q=e^{\frac{\pi i}{4}\frac{\tau-2}{\tau-1}}_{\mathstrut},
\quad\tau\to 1+i\,0\;\;\,\Big\}_{\mathstrut}\,,\\\\
\Big\{q=e^{\frac{\pi i}{4}\frac{\tau-1}{\tau +1}}_{\mathstrut},
\quad\tau\to -1+i\,0\Big\}\,.
\end{array}
$$
All the functions are holomorphic at  $\mathbb{H}^+$ and form the
field $\mathbb{C}(x,y)$ of meromorphic functions on the curve
$(1)$$:$
$R_1^{}\big(x(\tau)\big)+y(\tau)\,R_2^{}\big(x(\tau)\big)$.}

{\em Remark\;$4$\/}. Complete proof of the fact that all the series
are integer $q$-series goes beyond the scope of the present work. It
exploits some manipulations with Weierstrassian functions,
$\vartheta$-constants, and their transformations in
$\mathrm{PSL}_2(\mathbb{Z})$ \cite{br3}. As an example we exhibit
exact representations only for the series (\ref{XY}).

{\bf Proposition\;6.} {\em The integer series {\rm (\ref{XY})} have
an exact representation in form of products\/}$:$
\begin{equation}\label{XY2}
\left\{\!
\begin{array}{l}\ds
X=\qquad\,\big\{\!\pm\!1,\,\pm i \big\}\,\phantom{q}
\prod\limits_{k=1}^{\infty}\! \big( 1+q^{4k}\big)^2\big(1+q^{4k-2}
\big)^4
\\ \\\ds
Y=4\,\big\{\sqrt{\pm1},\,\sqrt{\pm i\,}\big\}\,q
\prod\limits_{k=1}^{\infty}\! \big(1+q^{4k} \big)^3\big(
1+q^{2k}\big)^6
\end{array}.
\right.
\end{equation}
These formulae can be considered as an alternative and much simpler
version of Burnside's parametrization itself: $Y^2=X^5-X$ is an
identity between the products (\ref{XY2}). Note that the series
(\ref{Xseries})--(\ref{Yseries}) are alternating ones and
exponential multipliers in the series (\ref{moon}) have been
introduced in order that the series be positive definite.

Since equation (\ref{X}) is a rational function of the quantities
$(q,\,X)$ and derivatives of $X$, there exist polynomial recurrence
relations for the coefficients of these canonical series. We do not
write down them here as they are rather cumbersome. See the work
\cite{rankin} for group properties of the $Y$-function and the work
\cite{br2} for complete $\vartheta$-treatment of Proposition\;6.

\subsection{Holomorphic Abelian integral as function of $\tau$}
Holomorphic Abelian integrals (as functions of $\tau$) are
fundamental objects in uniformization; though no really this fact
has been pointed out in the literature. All other Abelian integrals
and, in particular, meromorphic uniformizing functions, are
expressed through these holomorphic objects by means of Riemann
$\Theta$-functions. If  the uniformizing group is a Fuchsian group
of first kind \cite{ford} then all the integrals become
single-valued additive functions. The fact that the curve (1) covers
a torus suggests that the holomorphic integrals reduce to the
elliptic ones and therefore have representations in terms of
classical elliptic integrals.  In the case of Burnside's curve, as
it follows from  formulae (\ref{tor})--(\ref{cov}) and
(\ref{Palpha}), this leads, after some simplification, to the
following formula
\begin{equation}\label{holo}
\int\limits^{\,\,\,x}\!\frac{x\mp i\sqrt{i\,}}{\sqrt{x^5-x}}\,dx=
\sqrt{1+i\,}\,\wp^{-1}\!\!\left(
\mfrac{(1\pm\sqrt{2})\,(1-i)\,x}{\big(x-i\big)
\big(x+1\big)}-\mfrac{3\pm\sqrt{2}}{6};\mfrac53,\,
\mfrac{\mp7}{27}\sqrt{2} \right).
\end{equation}
On the other hand explicitly solvable Fuchsian equations on torus
described in Sec.\,IV.C allows one to construct $q$-series for these
integrals and to get other information.

{\bf Proposition\;7.} {\em The solution\/ $\alpha=\alpha(\tau)$ of
Schwarz's equation
$$
\begin{array}{l}
\ds\big[\alpha,\,\tau \big]=\ds -
\mfrac{1^{\mathstrut}}{2}\big(\wp(\alpha)+\wp(\alpha-\omega)+
\wp(\alpha-\omp)\big) -\mfrac{3}{8}
\big(\wp(\alpha-\Kappa)+\wp(\alpha+\Kappa)\big)+{}
\\\\
\ds \phantom{\ds\big[\alpha,\,\tau \big]=}
{}\,+\mfrac{9}{32}\sqrt[\leftroot{1}4]{8\,}\,i\,\big(
\zeta(\alpha-\Kappa)-\zeta(\alpha+\Kappa)\big)+ \mfrac{9}{16}
\big(\sqrt[\leftroot{1}4]{8\,}\,i\,\zeta(\Kappa)+2\sqrt{2}+2\big)\,,
\end{array}
$$
is a holomorphic additively automorphic function  with respect to
rank six group defined by monodromy of equation {\rm (\ref{new})}.
Arbitrary Abelian differential $d\alpha(\tau)$ is represented by the
two integer $q$-series\/}$:$
\begin{equation}\label{Adiff}
\left\{\!
\begin{array}{l}\ds
\phantom{X\,}\frac{dX}{Y}=\quad\;\;\;\big\{\sqrt{\pm1},\,\sqrt{\pm
i} \big\}\,2\prod\limits_{k=1}^{\infty}\! \big( 1-q^{4k}\big)
\big(1-q^{4k-2} \big)^2 \big( 1-q^{8k}\big)^3\cdot dq
\\ \\\ds
X\,\frac{dX}{Y}= \big\{\!\pm\!\sqrt{\pm1},\,-\sqrt{\mp i\,} \big\}\,
2 \prod\limits_{k=1}^{\infty}\! \big( 1-q^{4k}\big) \big(1+q^{4k-2}
\big)^2 \big( 1-q^{8k}\big)^3\cdot dq
\end{array}.
\right.
\end{equation}

{\em Proof}. By virtue of the formula (\ref{holo}) the quantity
$\alpha(\tau)$ is proportional to the holomorphic Abelian integral
on Burnside's Riemann surface and gives the inversion of the ratio
(\ref{Xi}) in the following explicit form
\begin{equation}\label{alpha}
\alpha^{\scriptscriptstyle\pm}(\tau)=\wp^{-1}\!\!\left(
\mfrac{(1\pm\sqrt{2})\,(1-i)\,x(\tau)}{\big(x(\tau)-i\big)
\big(x(\tau)+1\big)}-\mfrac{3\pm\sqrt{2}}{6};\mfrac53,\,
\mfrac{\mp7}{27}\sqrt{2} \right).
\end{equation}
Hence, being analytic function, it is everywhere finite in a domain
of its existence including the limiting points $\om_j$. The closed
loops surrounding the points $\alpha_j^{}=\{0,\om,\omp,\pm\Kappa\}$
determine five $2\!\times\!2$-matrix representations of generators
of automorphisms of the function $\alpha(\tau)$. These constitute a
subgroup. The shifts $\alpha\mapsto\alpha+\{2\,\om, 2\,\omp\}$ yield
two remaining generators of the full monodromy group. Rank of the
group (it is not free) is equal to six as there are a puncture and
one obvious relation between these seven elements of the monodromy.

Again, as in Proposition\;6, some routine manipulations with modular
$\vartheta$-forms and  (\ref{XY2}) lead to formulae (\ref{Adiff}).
Combining these formulae with (\ref{holo})--(\ref{alpha}) we get
$$
d\alpha^{\scriptscriptstyle\pm}(\tau)=\frac{X\mp
i\,\sqrt{i\,}}{\sqrt{1+i\,}}\,\frac{dX}{Y}\,.
$$
Series of a similar nature (without poles) can be obtained for
points $e_j^{}=\{0,\infty \}$.\hfill $\blacksquare$

We see that both the Puiseux series for multi-valued functions
considered in Sec.\,IV.B and Puiseux series for holomorphic
integrals (\ref{holo}) are transformed into single-valued series
when we use any of the global coordinates. In particular, we obtain
$q$-representation for complicated and seemingly chaotic Puiseux
series (\ref{puiseux}). The representation has a quite regular
structure defined by the two integer $q$-series (\ref{Adiff}):
$$
\begin{array}{l}
d\alpha^{\scriptscriptstyle+}(\tau)= M\,\big\{ 1- 2\,\gamma\,q^{2}-
q^{8}+ 6\,\gamma\,q^{10}- 6\,q^{16} -2\,\gamma\,q^{18}+
5\,q^{24}-4\,\gamma\,q^{26}+ 12\,q^{32}-{}
\\\\
\phantom{d\alpha^{\scriptscriptstyle+}(\tau)= 8\,i\,q\,\big\{ }
\!\!\!\!- 6\,q^{40}- 10\,\gamma\,q^{42}- 7\,q^{48}+
12\,\gamma\,q^{50}- 4\,q^{56}+6\,\gamma\,q^{58}+
\cdots\big\}\,dq\,,{}_{{}_{\ds\mathstrut}}
\end{array}
$$
\begin{equation}\label{ser}
\quad\alpha^{\scriptscriptstyle+}(\tau)=\omega+M\,\Big\{q-
2\,\gamma\,\mfrac{q^3}{3}- \mfrac{q^9}{9}+
6\,\gamma\,\mfrac{q^{11}}{11}-
6\,\mfrac{q^{17}}{17}-2\,\gamma\,\mfrac{q^{19}}{19}+
5\,\mfrac{q^{25}}{25}-4\,\gamma\,\mfrac{q^{27}}{27}+ \cdots
\Big\}\,,
\end{equation}
where
$$
\gamma= \sqrt{2}-1\,,\qquad M={\textstyle2\,\sqrt{\!\sqrt{2}+1}}
\,,\qquad q=e^{\frac{\pi i}{4}\frac{\tau-1}{2\tau
-1}}_{\mathstrut}\,.
$$
Note  that the term $(\cdots+ 12\,q^{32}- 6\,q^{40}- \cdots)$ is not
a misprint. The quantities $\gamma^1$ and $\gamma^0$ do not
alternate each other.

We remark that the complete group and Moonshine' treatment of this
nontrivial case of genus two is of independent interest. One
complication occurring in a higher-genera Moonshine is that there is
no canonical choice for Hauptmoduli \cite{gannon}. However the genus
two curves have a unique hyperelliptic shape and, as we have seen
now, Burnside's example can be thought of as the canonical in all
respects. What is more,  Moonshine' treatment, if any, of the {\em
holomorphic\/} ``toroidal Hauptmodulus'' $\alpha(\tau)$ (\ref{ser}),
its differential (\ref{Adiff}), and orbifold (\ref{new}) would be of
special interest.

Special attention must be given to the fact that once an explicit
formula  for uniformizing function, series  for holomorphic
integral, accessory parameters, or $\Psi$-function has been
obtained, all the series are recomputed one through  another. But
for the reasons pointed out above, the formula for holomorphic
integrals (in form of series (\ref{Adiff}), say) should be
considered as {\em the primary\/} object of the theory. This remains
valid even though we have no an explicit formula for the integrals
in terms of elliptic ones like (\ref{holo}); that is, when the cover
of torus does not exist.

\subsection{Modular forms for Burnside's function}
Poincar\'e's method of construction of automorphic functions as
ratios of  automorphic forms is widely known \cite{ford}. These
forms are analytic functions $\Theta(\tau)$ with the property:
$$
\Theta\!\left(\mfrac{a\,\tau+b}{c\,\tau+d}\right)=
(c\,\tau+d)^n\,\Theta(\tau)\,,
$$
where $n$ is the weight of the automorphic form $\Theta$ and
$(a,\,b,\,c,\,d)$ are substitutions of a group.  Burnside's example
fits in this classical construction and generates integer $q$-series
for modular forms.

{\bf Proposition\;8.} {\em The function of Burnside $x(\tau)$ is a
ratio of two automorphic holomorphic modular forms of weight $n=2$
with respect to  group $\boldsymbol{\Gamma}(4):$}
$$
x(\tau)=\frac{\Theta_1^{}(\tau)}{\Theta_2^{}(\tau)}\,,
\qquad\mbox{\rm where}\qquad \Theta_1^{}(\tau)=\wp(1)-\wp(2)\,,
\qquad \Theta_2^{}(\tau)=\wp(\tau)-\wp(2)\,.
$$

{\em Proof\/}. Derivative of any automorphic function is an
automorphic form of weight $n=2$ (Abelian differential) with respect
to automorphism group of the function. In our case, this is the
monodromy group of the equation (\ref{last}) --- the group
$\boldsymbol{\Gamma}(4)$ \cite{burnside,klein}. From (\ref{xt}) we
have that $\wp(2)$ is a form of weight $n=2$. From (\ref{p2}) we
have the same for $\wp(1)$ and, from (\ref{pt}), for $\wp(\tau)$.
Making use of the formulae (\ref{p2})--(\ref{pt}) we have
\begin{equation}\label{Theta1}
\Theta_1(\tau)= \frac94\,\frac{g_3^{}}{g_2^{}} \,
\frac{(x^3+x)(x^8+14\,x^4+1)}{x^{12}-33\,x^8-33\,x^4+1}=
\frac{1}{4}\,\frac{\pi\,i\,x_\tau}{x^2-1}\,.
\end{equation}
Correlation  the series/products (\ref{Xseries})--(\ref{XY2}) with
(\ref{Adiff})  (we omit details) shows that $\Theta_1(\tau)$ is
everywhere finite. The same is true for the form $\Theta_2(\tau)$.
\hfill $\blacksquare$

We exhibit some of such series only for the form $\Theta_1$. All the
representations can be derived from the formula (\ref{Theta1}). For
example the infinite point $\tau\to+i\,\infty$ with the coordinate
$q=e^{\frac{\pi i}{4}(\tau-3)}_{\mathstrut}$ yields the following
$q$-expansion:
$$
\begin{array}{l}
\ds \Theta_1(\tau)=-\frac{\pi^2}{16}\,
\frac{q\,x_{q}}{x^2-1}=\frac{\pi^2}{16}\!
\left.{}^{\ds\mathstrut}_{\ds\mathstrut} \right(
\!\!1+2\sum\limits_{k=1}^\infty q^{8k^2}\!\!
\left.{}^{\ds\mathstrut}_{\ds\mathstrut}
\right)^{\!\!4}={}\\\\
\ds \phantom{\Theta_1(\tau)}
=\frac{\pi^2}{16}\,\big(1+8\,q^8+24\,q^{16} +32\,q^{24} +
24\,q^{32}+48\,q^{40} +96\,q^{48}+64\,q^{56}+\cdots \big)\,.
\end{array}
$$
Zero of the form $\Theta_1$ corresponds to the cusp  (\ref{local})
and the following series
$$
\begin{array}{l}
\ds \Theta_1(\tau)=\frac{q\,x_{q}}{x^2-1}\,\ln^2\!q={}\\\\
\ds \phantom{\Theta_1(\tau)} = -4\,q^2\ln^2\!q\cdot\big(
1+4\,q^4+6\,q^8+8\,q^{12}+
13\,q^{16}+12\,q^{20}+14\,q^{24}+\cdots\big)\,.
\end{array}
$$
One of the explicit analytic expressions for this series is given by
the formula
$$
\Theta_1(\tau)=-4\,q^2\ln^2\!q\cdot\sum\limits_{k=0}^{\infty}\,
\sigma_1^{}(2\,k+1)\,q^{4k}\,,
$$
where $\sigma_1^{}(n)$ is sum of positive divisors of $n$. The form
$\Theta_2^{}(\tau)$ and all other points $\omega_j^{}$ are treated
in a similar manner but we do not write up here  exact
representations (which we have determined) to them in terms of
number-theoretic functions or Jacobi's $\vartheta$-constants. See
works \cite{mckay1,br2} for relevant information.

Alternatively, we could represent the function $x(\tau)$ in somewhat
unusual way. Namely in a form of ratio of two automorphic functions
$z=\frac{\wp_1^{}}{\wp_2^{}}-1,\;w=\frac{\wp_\tau^{}}{\wp_2^{}}-1$
rather than ratio of forms. Corresponding Schwarz's equations and
their monodromies for the functions $z$ and $w$ are derived from the
formulae (\ref{four})--(\ref{pt}). The following Schwarz's equation
elucidates completely this remark:
$$
z=2\,\frac{\wp_1^{}}{\wp_2^{}}+2 \qquad\Rightarrow\qquad
\big[z,\,\tau \big]=-\frac{27}{2}\,\frac{z^2+3}{z^2(z^2-9)^2}\,.
$$
Indeed, {\bf Aut}$\big(z(\tau)\big)=\boldsymbol{\Gamma}(2)\supset
\boldsymbol{\Gamma}(4)=\mbox{{\bf Aut}}\big(x(\tau)\big)$
\cite{burnside,klein,rankin,br2}.

\subsection{Conversion between uniformizations of different kinds}
The fact that both Whittaker's and Burnside's equations (\ref{last})
and (\ref{fuchs}) describe  uniformizations in their own rights
means that there is a one-to-one conformal transformation between
their global coordinates, i.\,e. {\em bi-holomorphic equivalence\/}.

Let $x(\tau)$ and $\boldsymbol{x}(\boldsymbol{\mu})$ be Burnside's
and Whittaker's uniformizing functions  respectively. We are looking
for a relation between $\tau$ and $\boldsymbol{\mu}$:
$\boldsymbol{\mu}=\boldsymbol\varphi(\tau)$. The function
$\boldsymbol\varphi$ and its inversion must be everywhere
holomorphic functions because we  deal with Fuchsian groups of first
kind: that is, groups having invariant circles \cite{ford,rankin2}
in the planes $(\tau)$ and $(\boldsymbol{\mu})$.

{\bf Proposition\;9.} {\em The coordinate $\boldsymbol{\mu}$ $($up
to a linear-fractional transformation$)$ is everywhere  holomorphic
function of\/ $\tau$ satisfying the following differential equation}
\begin{equation}\label{ff}
\boldsymbol\varphi:\quad\big[\tau,
\boldsymbol{\mu}\big]=-\frac{3}{8}\,\frac{g_2^{}(\tau)}{\pi^2}\,.
\end{equation}

{\em Proof\/}. We have $x(\tau)=\boldsymbol{x}(\boldsymbol{\mu})$
and
$$
\big[x,\tau \big]=-\frac12\,\frac{x^8+14\,x^4+1}{(x^5-x)^2}\,,
\qquad \big[\boldsymbol{x},\boldsymbol{\mu}
\big]=-\frac38\,\frac{\boldsymbol{x}^8+14\,\boldsymbol{x}^4+1}
{(\boldsymbol{x}^5-\boldsymbol{x})^2}\,.
$$
Invoking the transformation rule (\ref{prop}) of the object
$\boldsymbol{\mathfrak{D}}$ we deduce that
$$
\begin{array}{l}
\ds \big[x,\tau \big]=\big[\boldsymbol{x}(\boldsymbol{\mu}),\tau
\big]= \big[\boldsymbol{x}(\boldsymbol{\mu}),\boldsymbol{\mu}\big]+
\frac{1}{\boldsymbol{x}_{\boldsymbol{\mu}}^2}\,
\big[\boldsymbol{\mu},\tau \big]={}\\\\
\ds \phantom{\big[x,\tau \big]}=
-\frac38\,\frac{\boldsymbol{x}^8+14\,\boldsymbol{x}^4+1}
{(\boldsymbol{x}^5-\boldsymbol{x})^2}+
\frac{\boldsymbol{\mu}_\tau^2}{x_\tau^2}\,
\big[\boldsymbol{\mu},\tau \big]=
-\frac12\,\frac{x^8+14\,x^4+1}{(x^5-x)^2}\,.
\end{array}
$$
It follows that
$$
\big\{\boldsymbol{\mu},\tau \big\}=
-\frac18\,\frac{x^8+14\,x^4+1}{(x^5-x)^2}\,x_\tau{}^{\!\!\!\!2}\,.
$$
Correlating this expression with formulae (\ref{g2}) and (\ref{xt})
we arrive at formula (\ref{ff}).  The form $g_2^{}(\tau)$ is
everywhere finite  at $\mathbb{H}^+$ and hence the Schwarzian
$\{\boldsymbol{\mu},\tau\}$ and $\boldsymbol{\mu}$ itself are finite
as well. For the same reason function
$\tau=\boldsymbol{\varphi}^{-1}(\boldsymbol{\mu})$ has no fold
$a$-points (see (\ref{poles})) and is reversible into the function
of the same kind. \hfill $\blacksquare$

Explicitly solvable Schwarz's equations of the form
$[\tau,\boldsymbol{\mu}]=Q(\tau)$ with  holomorphic right-hand side
$Q(\tau)$ do exist. One of such examples is the nice equation
\begin{equation}\label{23}
\big[\tau,
\boldsymbol{\mu}\big]=-\frac23\,\frac{g_2^{}(\tau)}{\pi^2} \qquad
\Rightarrow\qquad
\frac{a\,\boldsymbol{\mu}+b}{c\,\boldsymbol{\mu}+d}=
\int\limits^{\,\,\tau}\! \widehat\eta^4(\tau)\,d\tau
\end{equation}
which comes from the following linear ODE with everywhere
holomorphic coefficients:
$$
\Psi_{\tau\mspace{-1mu}\tau}+
\frac{n+2}{\pi\,i}\,\eta(\tau)\,\Psi_\tau
-\frac{n}{6\pi^2}\,g_2^{}(\tau)\,\Psi=0\,.
$$
A direct check, with use of rules (\ref{diffg2g3}) and known
relation $\pi\ln_\tau\!\widehat\eta(\tau)=i\,\eta(\tau)$, shows that
$\Psi=\widehat\eta^{\,n\!}(\tau)$ solves this equation and  general
integral (\ref{23}) corresponds to the case $n=-2$.

A remarkable fact is that  equation (\ref{ff}) can also be exactly
integrated  but we shall present this material in a separate work.
Here we restrict ourselves to series representations. Analyzing this
equation one can show that in the neighborhood of infinity
$\tau\to+i\,\infty$ the global coordinate $q$ must have the form
$q=e^{\frac{\pi}{4}i\tau}_{\mathstrut}$. Invoking the well-known
expansion
$$
g_2^{}(\tau)=20\,\pi^4\!\left(
\mfrac{1}{240}+q^{8}+9\,q^{16}+28\,q^{24}+73\,q^{32}+
126\,q^{40}+252\,q^{48}+\cdots\right)
$$
we find that the series solutions have the following forms
$$
\begin{array}{rl}
\tilde{\boldsymbol{\mu}}\!\!\!&=
q-\frac{5}{21}\,q^9-\frac{78}{833}\,q^{17}+
\frac{4001}{39445}\,q^{25}+\frac{168948}{1711913}\,q^{33}+
\frac{42752022}{491319031}\,q^{41}- \cdots\,,\\\\
q\!\!\!&= \tilde{\boldsymbol{\mu}}+
\frac{5}{21}\,\tilde{\boldsymbol{\mu}}^9+\frac{503}{833}\,
\tilde{\boldsymbol{\mu}}^{17}+\frac{4138924}{2011695}\,
\tilde{\boldsymbol{\mu}}^{25}+
\frac{6383638315}{785768067}\,\tilde{\boldsymbol{\mu}}^{33}
+\cdots\,\qquad \big(\tilde{\boldsymbol{\mu}}=\frac{\pi
i}{4}\,\mu\big).
\end{array}
$$
Similar bi-holomorphic series exist for other cusps but  their
integer series realization, if any, is open question.

\section{Remarks concerning some of the literature}

Some of the expressions in the Secs.\,II, III and ground forms
appeared in (\ref{p2})--(\ref{g2})  have already occurred in the
literature \cite[p.\,59]{bolza}. For example $s^8-14\,s^4+1$, which
differs from $x^8+14\,x^4+1$ by a multiplier
$\sqrt[\uproot{2}\leftroot{1}4]{\!-1}$, appears throughout the
Schwarz Abhandlungen \cite{schwarz} in different contexts. Formula
(\ref{J}) was obtained by him \cite{schwarz}, Klein \cite{klein},
and Brioschi (1877) in relation to the groups of Platonic solids,
without mention of the explicit uniformization or function
$x(\tau)$. Slightly different (Legendre's) form of (\ref{J})
appeared in an earlier paper of Burnside on p.\,176 in {\em The
Messenger of Math. {\bf XXI}} $(1892)$. All this entails, in a
hidden form, some results of Secs.\,II and V. We should also remark
Schwarz's comments on pp.\,364--365 in \cite[{\bf II}]{schwarz},
where other candidates for a complete uniformization can be found.
We especially note the examples of Forsyth \cite[{\em Ex.}\,2--3,
p.\,188]{forsyth2}, \cite[{\em Ex.}\,13--14,
pp.\,242--243]{forsyth3} which  are not accompanied by any
references or comments however.

Some of the integer sequences presented above can be found in {\em
The On-Line Encyclopedia of Integer Sequences\/} by
N.\,J.\,A.\,Sloane which currently available at {\tt
http://www.research.att.com/\~{}njas/sequences}. Holomorphic and
meromorphic Abelian integrals associated with Whittaker's curve
(\ref{whit}) were considered by Legendre (with numerous details but
for the most part numerically) in {\em Troisi\`eme Suppl\'ement\/}
to \cite[see pp.\,207--269]{legendre}.

\section{Acknowledgements}
The author expresses a gratitude to Professor J.\,C.\,Eilbeck and
Professor E.\,Previato for much attention to the work and
hospitality in Heriot--Watt and Boston Universities where the most
part of the work was carried out.

The project was supported by a Royal Society/NATO Fellowship,
NSF/NATO-grant DGE--0209549, and Russian Science and Innovations
Federal Agency under contract 02.740.11.0238.



\thebibliography{99}

\bibitem{tah} \mbox{\sc
Aldrovandi, E. \& Takhtajan, L.\,A.} {\em Generating Functional in
CFT and Effective Action for Two-Dimensional Quantum Gravity on
Higher Genus Riemann Surfaces\/}. Comm. Math. Phys. {\bf188}, 29--67
(1997).

\bibitem{bolza} \mbox{\sc Bolza, O.}
{\em On binary sextics with linear transformations into
themselves\/}. Amer.\;Journ.\;Math.  {\bf 10}, 47--70 (1887).


\bibitem{burnside} \mbox{\sc Burnside, W.\,S.}
{\em Note on the equation \mbox{$y^2=x\,(x^4-1)$}\/}. Proc.\;London
Math.\;Soc.  {\bf XXIV}, 17--20  (1893).

\bibitem{br2} \mbox{\sc Brezhnev, Yu.\,V.}
\mbox{\em On the uniformization of algebraic curves\/}. Moscow Math.
Journ. {\bf 8}(2), 233--271 (2008).

\bibitem{br3} \mbox{\sc Brezhnev, Yu.\,V.}
\mbox{\em On functions of Jacobi and Weierstrass (I) and equation of
Painlev\'e\/}.\\ \texttt{http://arXiv.org/math.CA/0808.3486}.

\bibitem{chud2} \mbox{\sc Chudnovsky, D.\,V. \& Chudnovsky, G.\,V.}
{\em Transcendental Methods and Theta-Functions\/}.
Proc.\,Sympos.\,Pure Math. {\bf 49}(2), Part 2, 167--232  (1989).


\bibitem{chud} \mbox{\sc Chudnovsky, D.\,V. \& Chudnovsky, G.\,V.}
{\em Computer algebra in the service of mathematical physics and
number theory\/}. Computers in mathematics. Lect.\;Note in Pure and
Appl.\,Math. {\bf 125}, 109--232. Dekker: New York (1990).

\bibitem{dhar} \mbox{\sc Dhar, S.\,C.}{\em On the uniformization of a special kind
of algebraic curve of any genus\/}. Journ.\;London Math.\;Soc. {\bf
10}, 259--263 (1935).

\bibitem{erdei} \mbox{\sc Erd\'elyi, A., Magnus, W., Oberhettinger, F. \&
Tricomi, F.\,G.} {\em Higher transcendental functions\/}. {\bf 3}.
McGraw--Hill Book Comp., Inc., (1955).

\bibitem{ford} \mbox{\sc Ford, L.} {\em Automorphic Functions\/}.
McGraw--Hill Book Comp., Inc. New York (1929).

\bibitem{forsyth2} \mbox{\sc Forsyth, A.\,R.} {\em Theory of differential
equations\/}.  P.\,III ({\bf IV}): pp.\,1--534. Dover Publ., Inc.,
(1959).

\bibitem{forsyth3} \mbox{\sc Forsyth, A.\,R.} {\em A Treatise on differential
equations\/}. 6-th ed. Macmillan and Co., Ltd., (1943).

\bibitem{gannon} \mbox{\sc Gannon, T.}  {\em Monstrous Moonshine:
The first twenty-five years\/}. Bull.\,London Math.\,Soc. {\bf
38}(1), 1--33  (2006).


\bibitem{mckay2} \mbox{\sc Harnad, J. \& McKay, J.}
{\em Modular solutions to equations of generalized Halphen type\/}.
Proc. Royal Soc.\;London {\bf 456}(1994), 261--294  (2000).

\bibitem{jac} \mbox{\sc Jacobi, C.\,G.\,J.}
{\em Anzeige von Legendre th\'eorie des fonctions elliptiques\/}.
Crelle's Journal (1832), {\bf VIII}, 413--417; {\em Gesammelte
Werke\/} {\bf I}: 375--382. Berlin (1882--1891).

\bibitem{keen} \mbox{\sc Keen, L., Rauch, H.\,E. \& Vasquez, A.\,T.}
{\em Moduli of Punctured Tori and the Accessory Parameters of
Lam\'e's Equation\/}.  Trans.\,Amer.\,Math.\,Soc.  {\bf 255},
201--230 (1979).

\bibitem{klein} \mbox{\sc Klein, F. \& Fricke, R.} {\em Vorlesungen \"uber
die Theorie der elliptischen Modulfunktionen}. {\bf I}:
B.\,G.\,Teubner: Leipzig (1890).

\bibitem{kuusalo} \mbox{\sc Kuusalo, T. \& N\"a\"at\"anen, M.}
{\em Geometric uniformization in genus 2\/}.
Ann.\,Acad.\,Sci.\,Fenn. Ser. A. I. Math. {\bf 20}(2), 401--418
(1995).

\bibitem{legendre} \mbox{\sc Legendre, A.\,M.}
{\em Trait\'e des Fonctions Elliptiques et des Int\'egrales
Eul\'eriennes\/}. {\bf I--III}. Imprimere de Huzard-Courcier. Paris
(1825--1828).

\bibitem{eightfold} \mbox{\sc Levy, S. (Ed.)}
{\em The Eightfold Way. The Beauty of Klein's Quartic Curve.} MSRI
Publications, {\bf 35} (1998).

\bibitem{mckay1} \mbox{\sc McKay, J. \& Sebbar, A.}
{\em Fuchsian groups, automorphic functions and Schwarzians\/}.
Math.\;Annalen {\bf 318}(2), 255--275  (2000).

\bibitem{poincare} \mbox{\sc Poincar\'e, H.}
{\em Sur les groupes des \'equations lin\'eaires\/}. Acta
Mathematica {\bf 4}, 201--311  (1884).

\bibitem{rankin2} \mbox{\sc Rankin, R.\,A.}
{\em The Differential Equations Associated with the Uniformization
of Certain Alebraic Curves\/}. Proc.\,Royal Soc.\,Edinburgh. A {\bf
LXV}, 35--62 (1958).

\bibitem{rankin} \mbox{\sc Rankin, R.\,A.}
{\em Burnside's uniformization\/}. Acta Arithmetica {\bf 79}(1),
53--57  (1997).

\bibitem{schwarz} \mbox{\sc Schwarz, H.\,A.} {\em Gesammelte Mathematische
Abhandlungen\/}.  {\bf I--II}. Verlag von Julius Springer:  Berlin
(1890).

\bibitem{shabde} \mbox{\sc Shabde, N.\,G.}
{\em A note on automorphic functions\/}. Proc.\;Benares Math.\;Soc.
{\bf 16}, 39--46 (1934).

\bibitem{weber} \mbox{\sc Weber, H.}
{\em Ein Beitrag zu Poincar\'e's Theorie der Fuchs'schen
Functionen\/}. G\"ottinger Nachrichten,  359--370 (1886).

\bibitem{weyl} \mbox{\sc Weyl, H.}
{\em The concept of a Riemann surface\/}. 3-rd ed.
Addison--Wesley Publ.\;Comp., Inc. (1955).

\bibitem{whittaker1} \mbox{\sc Whittaker, E.\,T.}
{\em On the Connexion of Algebraic Functions with Automorphic
Functions\/}. Phil.\,Trans.\,Royal Soc.\,London  {\bf 192A}, 1--32
(1898).

\bibitem{whittaker2} \mbox{\sc Whittaker, E.\,T.}
{\em On hyperlemniscate functions, a family of automorphic
functions\/}. Proc.\;London Math.\;Soc. {\bf 4}, 274--278  (1929).

\bibitem{whittaker3} \mbox{\sc Whittaker, J.\,M.}
{\em The uniformisation of algebraic curves\/}. Proc.\;London
Math.\;Soc. {\bf 5}, 150--154  (1930).
\end{document}